\documentclass{aims} 
\usepackage{amsmath}
\usepackage{epsfig} 
\usepackage{epstopdf} 
\usepackage[colorlinks=true]{hyperref}
\usepackage{enumitem}
\setlist[itemize]{topsep=0pt,after=\vspace{1.5\baselineskip}}
\usepackage{array}
\usepackage{makecell}
\usepackage[utf8]{inputenc}     
\usepackage{ifthen} 
\usepackage{multicol}
\usepackage{multirow}
\usepackage{floatrow}
\usepackage{longtable}
\usepackage{caption}
\usepackage[utf8]{inputenc}     
\usepackage{ifthen} 
\usepackage{multicol}
\usepackage{multirow}
\usepackage{floatrow}
\usepackage{longtable}
\usepackage[singlelinecheck=false 
]{caption}
\usepackage{color, graphics}
\usepackage{graphicx, multicol, latexsym, amsmath, amssymb}
\usepackage{makecell}
\setcellgapes{4pt}
\usepackage{blindtext}
\usepackage{subfigure}
\usepackage{amsfonts}
\usepackage{enumitem}
\setlist[itemize]{topsep=0pt,after=\vspace{1.5\baselineskip}} 
\usepackage[colorlinks=true]{hyperref}
\usepackage{empheq}
\usepackage[T1]{fontenc}
\usepackage{tabularx}
\usepackage{pgfplots}
\usepackage{tikz}
\RequirePackage{doi}
\RequirePackage{filecontents}
\hypersetup{urlcolor=blue, citecolor=red}
\allowdisplaybreaks

\usepackage{color, graphics}
\usepackage{graphicx, multicol, latexsym, amsmath, amssymb}
\def\R{\mathbb R} \def\N{\mathbb N}

\def\R{\mathbb R} \def\N{\mathbb N} 
\def\TM{T_{\rm{max}}} 
\newtheorem{theorem}{Theorem}[section]

\newtheorem{lemma}[theorem]{Lemma}

\theoremstyle{definition}

\newtheorem{remark}{Remark}

{\begin{list}{}%
		{\setlength\labelwidth{0pt}%
			\setlength\leftmargin{#1}}\item[]}%
	{\end{list}}

\title[A nonlinear attraction-repulsion chemotaxis model with absorption]
{A nonlinear attraction-repulsion Keller--Segel model with double sublinear absorptions: criteria toward boundedness} 

\author[Yutaro Chiyo, Silvia Frassu and Giuseppe Viglialoro]{}


 \email{ycnewssz@gmail.com}
 \email{silvia.frassu@unica.it}
 \email{giuseppe.viglialoro@unica.it}

\thanks{YC is supported by \textit{JSPS KAKENHI}, Grant Number: JP22J11422. SF and GV are members of the Gruppo Nazionale per l'Analisi Matematica, la Probabilit\`a e le loro Applicazioni (GNAMPA) of the Istituto Nazionale di Alta Matematica (INdAM) and are partially supported by the research project \emph{Evolutive and Stationary Partial Differential Equations with a Focus on Biomathematics}, funded by Fondazione di Sardegna (2019). GV is partially supported by MIUR (Italian Ministry of Education, University and Research) Prin 2017 \emph{Nonlinear Differential Problems via Variational, Topological and Set-valued Methods} (Grant Number: 2017AYM8XW). }

\thanks{$^*$Corresponding author: giuseppe.viglialoro@unica.it}

\begin{document}
\maketitle

\centerline{\scshape Yutaro Chiyo}
\medskip
{\footnotesize
 \centerline{Department of Mathematics, 
    Tokyo University of Science
}
   \centerline{1-3, Kagurazaka, Shinjuku-ku, 
    Tokyo 162-8601}
   \centerline{Japan}
} 

\medskip

\centerline{\scshape Silvia Frassu and Giuseppe Viglialoro$^*$}
\medskip
{\footnotesize
 \centerline{Dipartimento di Matematica e Informatica}
   \centerline{Universit\`{a} di Cagliari, Via Ospedale 72, 09124. Cagliari}
   \centerline{Italy}
}

\bigskip



\begin{abstract}
This paper generalizes and extends to the case of nonlinear effects and logistic perturbations some results recently developed in \cite{FrassuGalvanViglialoro} where, for the linear counterpart and in absence of logistics, criteria toward boundedness for an attraction-repulsion Keller--Segel system with double saturation are derived. 
\end{abstract}

\tableofcontents
\section{Introduction and motivations}\label{Intro}
\subsection{General overview and state of the art}\label{IntroDiscussionStateOfTheArtSection}
The pioneering papers by Keller and Segel (\cite{K-S-1970,Keller-1971-MC,Keller-1971-TBC}), proposed in the 70's  to model the biological phenomena concerning with chemotaxis mechanisms, continue inspiring many researchers of the field toward the consideration of several related variants. 

In this regard, in the article we dedicate to a specific problem intimately connected to the forthcoming coupled system of Partial Differential Equations:
\begin{equation}\label{problemAttRep}
	 \hspace*{-0.2cm}
	\begin{cases}
		u_t= \nabla \cdot \left(D(u)\nabla u - S(u)\nabla v + T(u)
		\nabla w\right)+h(u)  & \text{ in } \Omega \times (0,\TM),\\
		\tau v_t= \Delta v -\beta  v +k(u,v)  & \text{ in } \Omega \times (0,\TM),\\
		\tau w_t= \Delta w-\delta w+l(u,w)& \text{ in } \Omega \times (0,\TM),\\
		u_{\nu}=v_{\nu}=w_{\nu}=0 & \text{ on } \partial \Omega \times (0,\TM),\\
		u(x,0)=u_0(x),\; \tau v(x,0)=\tau v_0(x),\; \tau w(x,0)=\tau w_0(x)& x \in \bar\Omega.\\ 
	\end{cases}
\end{equation}
Herein $\Omega$ is a bounded and smooth domain of $\R^n$, with
$n\in \N$, with boundary $\partial \Omega$, $\tau \in \{0,1\}$, $D=D(\xi), S=S(\xi), T=T(\xi), h=h(\xi), k=k(\xi,\eta)$ and
$l=l(\xi,\rho)$ are functions of their arguments $\xi,\eta\geq 0$ and $\rho\geq 0$ with a certain regularity and proper behavior. Moreover, further regular initial data  $u_0(x)\geq 0$, $ \tau v_0(x)\geq 0$ and $ \tau w_0(x)\geq 0$ are as well given, the subscript $\nu$ in $(\cdot)_\nu$ indicates the outward normal derivative on $\partial \Omega$, whereas $\TM \in (0,\infty]$ the maximal instant of time up to which solutions to the system do exist. 

In the frame of real biological phenomena, it is commonly used to indicate with $u=u(x,t)$, $v=v(x,t)$ and $w=w(x,t)$, respectively, the nonnegative value of a certain cell distribution (populations, organisms), of the concentration of a chemoattractant (i.e. a chemical signal whose effect is attracting the cells each other)  and of a chemorepellent density (i.e. a chemical signal with exactly the opposite effect of the chemoattractant); naturally the couple $(x,t)$ specifies the position and the instant of time where such values are considered. In this way, it is quite natural to interpret system \eqref{problemAttRep} as a model describing the motion of the cells, inside an insulated domain (zero-flux on the border: homogeneous Neumann boundary conditions), under the flux effect $D(u)\nabla u - S(u)\nabla v + T(u)\nabla w$, combination of diffusive ($D(u)$), attractive ($S(u)$) and repulsive ($T(u))$ impacts, and the external action of a source $h(u)$. As expected, the effect of such impacts are intimately connected to the expression of the diffusion, the attraction and the repulsion, while the source may increment, decrement or both the cells' density. Further, the attractive and repulsive signals growths evolve conforming the rates $k(u,v)$ and $l(u,w)$, respectively (second and third equation in \eqref{problemAttRep}). 
\begin{remark}
Even though the analysis of this paper is merely theoretical, it is worthwhile emphasizing that model \eqref{problemAttRep} finds some applications in the inflammation observed in Alzheimer's disease when $h(u)\equiv 0$, $k(u,v)$ and $l(u,w)$ are linear functions only of the cell density (evolutive PDEs for chemical in this form are also known as signal-production models). More precisely, \cite{Luca2003Alzheimer} deals with the description of the gathering mechanisms of microglia and dimensional, numerical and experimental analyses, in bounded intervals are proposed.
\end{remark}
As expected, the cellular motility is extremely sensitive to the actions of the above factors governing model \eqref{problemAttRep}: in particular, in general no stable behavior is conceivable but, conversely, even weak changes in the related values may strongly influence the dynamics. In the specific, the evolution might relax toward global stabilization and convergence to equilibrium of the cell distribution $u$, but it could even degenerate into the so-called \textit{chemotactic collapse}, the mechanism resulting in uncontrolled aggregation processes for $u$, eventually blowing up/exploding at finite time. By the mathematical point of view, in the first case solutions $(u,v,w)$ are defined and bounded for all $(x,t)$ in $\Omega \times (0,\infty)$, in the other for a certain finite time $\TM$, the solution $(u,v,w)$ becomes unbounded approaching $\TM$.

Let us mention now some known results in the literature dealing with signal-production models coming from \eqref{problemAttRep}; in this situation, $k$ and $l$ are positive functions only of $u$, and as $u$ itself increases so $v$ and $w$ do: 
\begin{enumerate}
	\item For $\tau=0$ in the equations for $v$ and $w$, $h(u)\equiv 0$ and for the  linear case (i.e., diffusion $D(u)=u$, chemosensitivities $S(u)=T(u)=u$ and chemoattractant and  chemorepellent $k(u,v)=\alpha u$, $\alpha>0$, and $l(u,w)=\gamma u$, $\gamma>0$, respectively), the value $\xi\gamma-\chi\alpha$ is critical for $n=2$: particularly, if $\xi\gamma-\chi\alpha>0$	(repulsion dominates on attraction), in any dimension all solutions to the model are globally bounded, whereas for
	$\xi\gamma-\chi\alpha<0$ (attraction dominates repulsion) unbounded solutions can be detected (see
	\cite{GuoJiangZhengAttr-Rep,LI-LiAttrRepuls,TaoWanM3ASAttrRep,VIGLIALORO-JMAA-BlowUp-Attr-Rep}). Indeed, when $\xi\gamma-\chi\alpha>0$ and for $h(u)=k u-\mu u^\beta$, $k\in \R,\mu>0$ and some $\beta>1$, in \cite{ChiyoMarrasTanakaYokota2021} some questions on the blow-up scenario are discussed.  
	\item For more general production laws, respectively $k(u,v)\approx \alpha u^s$, $s>0,$ and $l(u,w)\approx\gamma u^r$, $r\geq 1$ (whereas $D,S, T$ and $\tau$ as before) the following is available in the literature  (\cite{ViglialoroMatNacAttr-Repul}): for $n\geq 2$, for every  $\alpha,\beta,\gamma,\delta,\chi>0$,  and for $r>s\geq 1$ (resp. $s>r\geq 1$), there exists $\xi^*>0$ (resp. $\xi_*>0$) such that if $\xi>\xi^*$ (resp. $\xi\geq \xi_*$), any sufficiently regular initial distribution $u_0(x)\geq 0$ (resp. $u_0(x)\geq 0$ enjoying some smallness assumptions) infers a unique classical and bounded solution. In addition the same conclusion holds true for every  $\alpha,\beta,\gamma,\delta,\chi,\xi>0$, $0<s<1$, $r=1$ and  any sufficiently regular $u_0(x)\geq0$. 
\end{enumerate}
The chemotactic collapse appearing in the signal-production situations, apparently (but this is an open question in the field yet), cannot occur for alike signal-absorption models, where to high values of $u$ correspond small ones (eventually vanishing) of $v$ and $w$. As a matter of fact, even for the so called Keller--Segel system with consumption, which is the simplified two-unknowns version of problem \eqref{problemAttRep} reading as
\begin{equation}\label{problemOriginalKSCosnumption}
	\begin{cases}
		u_t= \Delta u - \chi \nabla \cdot (u \nabla v)  & \text{ in } \Omega \times (0,\TM),\\
		v_t=\Delta v-u v & \text{ in } \Omega \times (0,\TM),\\
		u_{\nu}=v_{\nu}=0 & \text{ on } \partial \Omega \times (0,\TM),\\
		u(x,0)=u_0(x), \; v(x,0)= v_0(x) & \;x \in \bar\Omega,
	\end{cases}
\end{equation}
the occurrence of blow-up has not been found. In particular all 
classical solutions $(u,v)$ to \eqref{problemOriginalKSCosnumption}
are  uniformly bounded in one and two-dimensional
settings, independently by some size of $(u_0,v_0)$:  the case $n=1$ can be justified by standard procedures, while for $n=2$ the result is a consequence of more general analyses discussed in \cite{WinklerN-Sto_CPDE}
and \cite{WinklerN-Sto_2d}. Conversely, for $n\geq 3$, boundedness requires the
smallness assumption $\chi \lVert v_0\lVert_{L^\infty(\Omega)}\leq
\frac{1}{6(n+1)}$, as proved in \cite{TaoBoun}. (This condition is improved in  $\chi \lVert v_0\lVert_{L^\infty(\Omega)}<
\frac{\pi}{\sqrt{2(n+1)}}$; see \cite{BaghaeiKhelghatiIranian}.) On the other hand, if in models like \eqref{problemOriginalKSCosnumption} some smoothing effects on the dynamics of problem
 are introduced, boundedness of related solutions remains valid even when $\chi \lVert
v_0\lVert_{L^\infty(\Omega)}$ is larger than the mentioned values. In the specific: 
\begin{enumerate}[resume]
	\item In \cite{LankeitWangConsumptLogistic}, where the motion of the cells is affected by a logistic source, so reading
	\begin{equation*}
		u_t= \Delta u - \chi \nabla \cdot (u \nabla v) +ku-\mu u^2, \quad
		 k,\mu>0,
	\end{equation*}
	the authors establish that the resulting Cauchy problem admits classical bo\-un\-ded solutions for
	arbitrarily large  $\chi \lVert v_0\lVert_{L^\infty(\Omega)}$ provided $\mu$ is also larger than a certain expression  increasing with $\chi \lVert v_0\lVert_{L^\infty(\Omega)}$;
	\item In \cite{MarrasViglialoroMathNach}, together with the dampening logistic term the equation for the particles' density is also perturbed by nonlinear diffusion and sensitivity:
	\begin{equation*}
		\hspace{1cm}
		u_t= \nabla \cdot \left((u+1)^{m_1-1}\nabla u -u (u+1)^{m_2-1}\nabla v\right) +ku-\mu u^2,\quad m_1,m_2 \in \R. 
	\end{equation*}
	Similarly to the previous case, globality and boundedness are derived whenever the smoothness parameter $\mu$ of the logistic is large enough and the diffusion dominates the attraction, in the sense that $m_1>2m_2-1$.
	\end{enumerate}
Naturally, in the frame of chemotaxis models with two signals, beyond the double-production cases aforementioned, double-saturation or consumption-production mechanisms can be considered. In this sense, and always with reference to   \eqref{problemAttRep}:  
\begin{enumerate}[resume]
		\item When $k(u,v)\approx \beta v-u^\alpha v$ in the second equation (with $\tau=1$) and  $l(u,w)\approx u^l$ in the third (with $\tau=0$) are fixed, the cells' density produces chemorepellent and absorbs chemoattractant;  in \cite{frassuviglialoroConsumptionProduction} boundedness is established (i) for $l=1$,
	$n\in \{1,2\}$, $\alpha\in (0,\frac{1}{2}+\frac{1}{n})\cap (0,1)$ and any $\xi>0$, (ii) for
	$l=1$, $n\geq 3$, $\alpha\in (0,\frac{1}{2}+\frac{1}{n})$ and $\xi$ larger than a quantity
	depending on $\chi \lVert v_0 \rVert_{L^\infty(\Omega)}$, (iii) for $l>1$ any $\xi>0$, and
	in any dimension; 
		\item \label{ItemFrassuGalvanViglialoro} For $\tau=1$ in the equations for $v$ and $w$, $h(u)\equiv 0$ and for the  linear case (i.e., diffusion $D(u)=u$, chemosensitivities $S(u)=\chi u$ and $T(u)=\xi u$, $\chi,\xi>0$) and chemoattractant and chemorepellent $k(u,v)=\beta v - uv, l(u,w)=\delta w - uw$ (double-signal saturation), for $n\geq 3$ all solutions emanating from sufficiently regular data such that  $0<\chi \lVert v_0 \rVert_{L^\infty(\Omega)} <\frac{1}{5n}$ and $0<\xi \lVert w_0 \rVert_{L^\infty(\Omega)} <\frac{1}{5n}$ are globally bounded. This is proved in \cite{FrassuGalvanViglialoro} (where also two- and three-dimensional numerical simulations are presented) and through the lines of this paper we will give some more hints on this research since the present article represents a generalization of what derived in \cite{FrassuGalvanViglialoro}. 
		\item \label{ItemFrassuLiViglialoro} For $\tau=1$ in the equation for $v$ and $k(u,v)=\beta v-uv$ (chemoattractant consumed), for $\tau=0$ in the third and $l(u,w)\approx u^l$, $l\geq 1$ (chemorepellent produced), for the diffusion $D(u)\approx u^{m_1}$, chemosensitivities $S(u)\approx u^{m_2}$ and $T(u)\approx u^{m_3}$ ($m_1,m_2,m_3\in \R$), for $h(u)\approx ku -\mu u^\beta$ ($k\in \R$, $\mu>0$ and $\beta>1$), for $n\geq 3$ all solutions emanating from sufficiently regular data $u_0,v_0$ are globally bounded provided $m_1>\varphi(m_2,m_3,n,\alpha,\beta,l)$, for specific expressions of $\varphi$. 
	\end{enumerate}
\section{Main claims and organization of the paper}\label{ClamimsSection}
\subsection{The model: presentation and some notations}\label{IntroDiscussionSection}
Let $\Omega \subset \R^n$, $n \geq 2$, be a bounded and smooth domain, $\chi,\xi>0$, $m_1,m_2,m_3\in\R$, $f(u), g(u)$ and $h(u)$ be reasonably regular functions generalizing the prototypes $f(u)=K_1 u^\alpha$, $g(u)=K_2 u^\gamma$, and $h(u)=k u - \mu u^{\beta}$ with $K_1, K_2, \mu>0$, $k \in \R$ and suitable $\alpha, \gamma, \beta>0$. Once nonnegative initial configurations $u_0$, $v_0$ and $w_0$ are fixed, herein we are interested in this nonlinear attraction-repulsion chemotaxis model, naturally obtainable as a particular case of problem \eqref{problemAttRep}: 
\begin{equation}\label{problem}
	\begin{cases}
		u_t= \nabla \cdot ((u+1)^{m_1-1}\nabla u - \chi u(u+1)^{m_2-1}\nabla v\\
		\hspace{7mm} {}+\xi u(u+1)^{m_3-1}\nabla w) + h(u) & \textrm{in } \Omega \times (0,T_{\rm{max}}),\\
		v_t=\Delta v-f(u)v  & \textrm{in } \Omega \times (0,T_{\rm{max}}),\\
		w_t= \Delta w - g(u)w & \textrm{in } \Omega \times (0,T_{\rm{max}}),\\
		u_{\nu}=v_{\nu}=w_{\nu}=0 & \textrm{on } \partial \Omega \times (0,T_{\rm{max}}),\\
		u(x,0)=u_0(x), \; v(x,0)=v_0(x), \; w(x,0)=w_0(x) & x \in \bar\Omega.
	\end{cases}
\end{equation}
In the light of what previously said, it should be easy to convince ourselves that this model  brings together many of the characteristics above discussed: nonlinear diffusion, sensitivities and growth rates, as well as general logistic terms. 

To our knowledge, the literature provides partial results connected to model \eqref{problem} only for the case $m_1=m_2=m_3=1$,  $f(u)=g(u)=u$ and $h(u)=0$ (see \cite{FrassuGalvanViglialoro}); the attained results were summarized in item \ref{ItemFrassuGalvanViglialoro}. In particular, system \eqref{problem} appears as a natural continuation of model described in item \ref{ItemFrassuGalvanViglialoro} itself, and it is worthwhile developing a general $n$-dimensional analysis in order to extend the mathematical comprehension. Specifically, we aim at deriving sufficient conditions involving the parameters used in problem \eqref{problem} according to which it admits classical solutions which are global and uniformly bounded in time. Specifically, we look for nonnegative functions $u=u(x,t), v=v(x,t), w=w(x,t)$ defined for $(x,t) \in \bar{\Omega}\times [0,T_{\rm{max}})$, and $T_{\rm{max}}=\infty$, such that 
\begin{align}\label{ClassicalAndGlobability}
	\begin{split}
		u, v, w\in C^0(\bar{\Omega}\times [0,\infty))\cap  C^{2,1}(\bar{\Omega}\times (0,\infty)) 
		\cap L^\infty((0, \infty);L^{\infty}(\Omega)),
	\end{split}
\end{align}
and satisfying for all $(x,t)\in \bar{\Omega}\times[0,\infty)$  all the relations in  \eqref{problem}.

To this aim, we require that $f$, $g$ and $h$ comply with 
\begin{align}\label{f}
	\begin{split}
		f,g \in C^1(\R) \quad \textrm{with} \quad   &0\leq f(s)\leq K_1 s^{\alpha}  \textrm{ and } 0\leq g(s)\leq K_2 s^{\gamma}, \\
		&\textrm{for some} \quad  K_1, K_2, \alpha, \gamma >0 \quad \textrm{and all } s \geq 0,
	\end{split}
\end{align}
and 
\begin{align}\label{h}
	\begin{split}
		h \in C^1(\R) \quad \textrm{with} \quad &h(0)\geq 0  \textrm{ and } h(s)\leq k s-\mu s^{\beta}, \\ 
		&\textrm{for some}\quad k \in \R,\,\mu>0,\, \beta>1\, \quad \textrm{and all } s \geq 0.
	\end{split}
\end{align}

\begin{remark}
	As the reader can expect, the results below will depend on the parameters $\alpha, \gamma, m_1, m_2, m_3$ and $n$. In particular, in view of the formulations of the second and third equation (which are somehow  ``exchangeable''), by permuting some of those constants a number of assumptions connecting their values can be seen as symmetric case of other. In this sense, even though the forthcoming presentation may appear hard to read, we want to underline that we made an important effort to include all the aforementioned permutations in the clearest way. This is the reason why we dedicate a part of the manuscript to define crucial constants in the computation; this is precisely what $\S$\ref{NotazioniSection} includes.   
\end{remark}

\subsubsection{Notations}\label{NotazioniSection}
We will make reference to these quantities in the absence of the logistic term:
\begin{equation*}
	\mathcal{A}:= \min
	\left\{\!\begin{aligned}
		&\max\Bigl\{2m_2-1, 2m_3-1, \frac{n-2}{n}\Bigr\}, \,\max\Bigl\{m_2-\frac{1}{n}, m_3-\frac{1}{n}, \frac{n-2}{n}\Bigr\}\\[1ex] 
		&\max\Bigl\{2m_2-1, m_3-\frac{1}{n}, \frac{n-2}{n}\Bigr\},\,
		\max\Bigl\{m_2-\frac{1}{n}, 2m_3-1, \frac{n-2}{n}\Bigr\}\\[1ex]
		&\quad \quad \quad \quad \quad \quad \quad \quad \quad  \max\Bigl\{m_2-\frac{1}{n}, m_3-\frac{1}{n}\Bigr\}
	\end{aligned}\right\}, 
\end{equation*}

\begin{equation*}
	\mathcal{B} := \min
	\left\{\!\begin{aligned}
		&\max\Bigl\{m_2-\frac{2}{n}+\alpha, m_3-\frac{2}{n}+\gamma\Bigr\},\,\max\Bigl\{2m_2, 2m_3, \frac{n-2}{n}\Bigr\}\\[1ex]
		&\max\Bigl\{m_2-\frac{2}{n}+\alpha, 2m_3, \frac{n-2}{n}\Bigr\},\,\max\Bigl\{2m_2, m_3-\frac{2}{n}+\gamma,\frac{n-2}{n}\Bigr\}
	\end{aligned}\right\},
\end{equation*}
\begin{equation*}
	\mathcal{C} := \min
	\left\{\!\begin{aligned}
		&\max\Bigl\{m_2+\frac{n\alpha-2}{n\alpha-1}, m_3 + \frac{n\gamma-2}{n\gamma-1}\Bigr\},\,\max\Bigl\{2m_2, 2m_3, \frac{n-2}{n}\Bigr\}\\[1ex]
		&\max\Bigl\{m_2+\frac{n\alpha-2}{n\alpha-1}, 2m_3, \frac{n-2}{n}\Bigr\},\,\max\Bigl\{2m_2,m_3 + \frac{n\gamma-2}{n\gamma-1}, \frac{n-2}{n}\Bigr\}
	\end{aligned}\right\}, 
\end{equation*}
\begin{equation*}
	\mathcal{D}:= \min
	\left\{
		\max\Bigl\{m_2+\frac{n\alpha-2}{n\alpha-1}, m_3 + \frac{n\gamma-2}{n\gamma-1}\Bigr\},\,\max\Bigl\{m_2+\frac{n\alpha-2}{n\alpha-1}, 2m_3, \frac{n-2}{n}\Bigr\}
\right\},
\end{equation*}
\begin{equation*}
	\mathcal{E} := \min
	\left\{
		\max\Bigl\{m_2+\frac{n\alpha-2}{n\alpha-1}, m_3 + \frac{n\gamma-2}{n\gamma-1}\Bigr\},\,\max\Bigl\{2m_2, \frac{n-2}{n}, m_3+\frac{n\gamma-2}{n\gamma-1}\Bigr\}
\right\}, 
\end{equation*}
\begin{equation*}
	\mathcal{F} := \max\Bigl\{m_2+\frac{n\alpha-2}{n\alpha-1}, m_3 + \frac{n\gamma-2}{n\gamma-1}\Bigr\}.
\end{equation*}
\begin{equation*}
	\mathcal{G}:= \min
	\left\{\!\begin{aligned}
		&\max\Bigl\{m_2-\frac{1}{n}, m_3-\frac{2}{n}+\gamma\Bigr\},\,\max\Bigl\{2m_2-1, 2m_3, \frac{n-2}{n}\Bigr\}\\[1ex]
		&\max\Bigl\{m_2-\frac{1}{n}, 2m_3, \frac{n-2}{n}\Bigr\},\,\max\Bigl\{2m_2-1,m_3-\frac{2}{n}+\gamma,\frac{n-2}{n}\Bigr\}\\[1ex]
		&\quad \quad \quad \quad  \quad \quad  \max\Bigl\{m_2-\frac{1}{n}, m_3-\frac{2}{n}+\gamma,\frac{n-2}{n}\Bigr\}
	\end{aligned}\right\},
\end{equation*}
\begin{equation*}
	\mathcal{H}:= \min
	\left\{\!\begin{aligned}
		&\max\Bigl\{m_2-\frac{1}{n}, m_3 + \frac{n\gamma-2}{n\gamma-1}\Bigr\},\,\max\Bigl\{2m_2-1, 2m_3, \frac{n-2}{n}\Bigr\}\\[1ex]
		&\max\Bigl\{m_2-\frac{1}{n}, 2m_3, \frac{n-2}{n}\Bigr\},\,\max\Bigl\{2m_2-1,m_3 + \frac{n\gamma-2}{n\gamma-1},\frac{n-2}{n}\Bigr\}\\[1ex]
		&\quad \quad \quad  \quad \quad  \quad \max\Bigl\{m_2-\frac{1}{n}, m_3 + \frac{n\gamma-2}{n\gamma-1}, \frac{n-2}{n}\Bigr\}
	\end{aligned}\right\},
\end{equation*}
\begin{equation*}
	\mathcal{I}:= \min
	\left\{\!\begin{aligned}
		&\max\Bigl\{m_2-\frac{1}{n}, m_3 + \frac{n\gamma-2}{n\gamma-1}\Bigl\},\,\max\Bigl\{2m_2-1, m_3 + \frac{n\gamma-2}{n\gamma-1},\frac{n-2}{n}\Bigl\}\\[1ex]
		&\quad \quad \quad \quad \quad \quad \quad \max\Bigl\{m_2-\frac{1}{n}, m_3 + \frac{n\gamma-2}{n\gamma-1},\frac{n-2}{n}\Bigl\}
	\end{aligned}\right\},
\end{equation*}
\begin{equation*}
	\mathcal{J}:= \min
	\left\{\!\begin{aligned}
		&\max\Bigl\{m_2-\frac{2}{n}+\alpha, m_3 + \frac{n\gamma-2}{n\gamma-1}\Bigr\},\,\max\Bigl\{2m_2, 2m_3, \frac{n-2}{n}\Bigr\}\\[1ex]
		&\max\Bigl\{m_2-\frac{2}{n}+\alpha, 2m_3, \frac{n-2}{n}\Bigr\},\,\max\Bigl\{2m_2, m_3 + \frac{n\gamma-2}{n\gamma-1},\frac{n-2}{n}\Bigr\}
	\end{aligned}\right\},
\end{equation*}
\begin{equation*}
	\mathcal{K} := \min
	\left\{
		\max\Bigl\{m_2-\frac{2}{n}+\alpha, m_3 + \frac{n\gamma-2}{n\gamma-1}\Bigr\},\,\max\Bigl\{2m_2,\frac{n-2}{n}, m_3 + \frac{n\gamma-2}{n\gamma-1}\Bigr\}
\right\}.
\end{equation*}
For the logistic case we will refer to these quantities:
\begin{equation*}
	\mathcal{A'}:= \min
	\left\{\!\begin{aligned}
		&\max\Bigl\{2m_2-1, 2m_3-1, \frac{n-2}{n}\Bigr\},\,\max\Bigl\{m_2-\frac{1}{n}, m_3-\frac{1}{n}, \frac{n-2}{n}\Bigr\}\\[1ex]
		&\max\Bigl\{2m_2-1, m_3-\frac{1}{n}, \frac{n-2}{n}\Bigr\},\,\max\Bigl\{m_2-\frac{1}{n}, 2m_3-1, \frac{n-2}{n}\Bigr\}\\[1ex]
		&\max\Bigl\{2m_2-\beta, 2m_3-\beta, \frac{n-2}{n}\Bigr\},\,\max\Bigl\{2m_2-\beta, 2m_3-\beta\Bigr\}\\[1ex]
		&\max\Bigl\{m_2-\frac{1}{n}, 2m_3-\beta,\frac{n-2}{n}\Bigr\},\,\max\Bigl\{2m_2-1,2m_3-\beta,\frac{n-2}{n}\Bigr\}\\[1ex]
		&\max\Bigl\{2m_2-\beta, 2m_3-1, \frac{n-2}{n}\Bigr\},\,\max\Bigl\{2m_2-\beta, m_3-\frac{1}{n}, \frac{n-2}{n}\Bigr\}
	\end{aligned}\right\}, 
\end{equation*}
\begin{equation*}
	\mathcal{B'} := \min
	\left\{\!\begin{aligned}
		&\max\Bigl\{2m_2, 2m_3, \frac{n-2}{n}\Bigr\},\,
		\max\Bigl\{2m_2+1-\beta, 2m_3+1-\beta\Bigr\}\\[1ex]
		&\max\Bigl\{2m_2, 2m_3+1-\beta,\frac{n-2}{n}\Bigr\},\,
		\max\Bigl\{2m_2+1-\beta, 2m_3,\frac{n-2}{n}\Bigr\}
	\end{aligned}\right\},
\end{equation*}
\begin{equation*}
	\mathcal{C'} := \min
	\left\{\!\begin{aligned}
		&\max\Bigl\{2m_2-1, \frac{n-2}{n}, 2m_3\Bigr\},\,\max\Bigl\{2m_2-1, \frac{n-2}{n}, 2m_3+1-\beta\Bigr\}\\[1ex]
		&\max\Bigl\{m_2-\frac{1}{n}, 2m_3,\frac{n-2}{n}\Bigr\},\,\max\Bigl\{m_2-\frac{1}{n},\frac{n-2}{n},2m_3+1-\beta\Bigr\}\\[1ex]
		&\max\Bigl\{2m_2-\beta, \frac{n-2}{n}, 2m_3\Bigr\},\,\max\Bigl\{2m_2-\beta, \frac{n-2}{n}, 2m_3+1-\beta\Bigr\}\\[1ex]
		&\quad \quad \quad \quad  \quad \quad \quad \quad \max\Bigl\{2m_2-\beta, 2m_3+1-\beta\Bigr\}
	\end{aligned}\right\}.
\end{equation*}
\subsection{Statements of the theorems and discussions}
With reference to the notations introduced in $\S$\ref{NotazioniSection}, let us now give the conclusions proved in this paper.
\begin{theorem}[The non-logistic case]\label{MainTheorem}
	Let $\Omega$ be a smooth and bounded domain of $\mathbb{R}^n$, with $n\geq 2$, $\chi, \xi$ positive reals, and $h \equiv 0$. Moreover, for $\alpha, \gamma >0$ and $m_1, m_2, m_3 \in \R$, let $f$ and $g$ fulfill \eqref{f} for each of the following cases:

	\begin{minipage}[t]{0.5\textwidth}
		\begin{enumerate}[label=$A_{\arabic*}$)]
			\setlength{\itemsep}{1.5mm}
			\item \label{A1} $\alpha, \gamma \in \left(0, \frac{1}{n}\right]$, $m_1>\mathcal{A}$;
			\item \label{A2} $\alpha, \gamma \in \left(\frac{1}{n},\frac{2}{n}\right)$,  $m_1>\mathcal{B}$;
			\item \label{A3} $\alpha, \gamma \in \left[\frac{2}{n},1\right)$, $m_1>\mathcal{C}$;
			\item \label{A4} $\alpha \in \left[\frac{2}{n},1\right], \gamma \in \left[\frac{2}{n},1\right)$, $m_1>\mathcal{D}$;
			\item \label{A5} $\alpha \in \left[\frac{2}{n},1\right), \gamma \in \left[\frac{2}{n},1\right]$,  $m_1>\mathcal{E}$;
			\item \label{A6} $\alpha, \gamma \in \left[\frac{2}{n},1\right]$, $m_1>\mathcal{F}$;
			\item \label{A7} $\alpha \in \left(0,\frac{1}{n}\right], \gamma \in \left(\frac{1}{n},\frac{2}{n}\right)$, $m_1>\mathcal{G}$;
			\item \label{A8} $\alpha \in \left(0,\frac{1}{n}\right], \gamma \in \left[\frac{2}{n},1\right)$; 
			$m_1>\mathcal{H}$;
		\end{enumerate}
	\end{minipage}
	\hspace{-4mm}
	\begin{minipage}[t]{0.5\textwidth}
		\begin{enumerate}[label=$A_{\arabic*}$)]
			\setlength{\itemsep}{1.5mm}
			\setcounter{enumi}{8}
			\item \label{A9} $\alpha \in \left(0,\frac{1}{n}\right], \gamma \in \left[\frac{2}{n},1\right]$, $m_1>\mathcal{I}$;
			\item \label{A10} $\alpha \in \left(\frac{1}{n},\frac{2}{n}\right), \gamma \in \left(0,\frac{1}{n}\right]$, $m_1>\mathcal{G}^t$;
			\item \label{A11} $\alpha \in \left(\frac{1}{n},\frac{2}{n}\right), \gamma \in \left[\frac{2}{n},1\right)$, $m_1>\mathcal{J}$;
			\item \label{A12} $\alpha \in \left(\frac{1}{n},\frac{2}{n}\right), \gamma \in \left[\frac{2}{n},1\right]$, $m_1>\mathcal{K}$;
			\item \label{A13} $\alpha \in \left[\frac{2}{n},1\right), \; \gamma \in \left(0,\frac{1}{n}\right]$, $m_1>\mathcal{H}^t$;
			\item \label{A14} $\alpha \in \left[\frac{2}{n},1\right], \; \gamma \in \left(0,\frac{1}{n}\right]$, $m_1>\mathcal{I}^t$;
			\item \label{A15} $\alpha \in \left[\frac{2}{n},1\right), \; \gamma \in \left(\frac{1}{n},\frac{2}{n}\right)$,  $m_1>\mathcal{J}^t$;
			\item \label{A16} $\alpha \in \left[\frac{2}{n},1\right], \; \gamma \in \left(\frac{1}{n},\frac{2}{n}\right)$, $m_1>\mathcal{K}^t$, 
		\end{enumerate}
	\end{minipage}

	\medskip
	\noindent
	the superscript $t$ denoting the case where the roles of $\alpha$ and $m_2$ are taken by $\gamma$ and $m_3$, respectively. Then for any initial data $(u_0,v_0,w_0)\in (W^{1,\infty}(\Omega))^3$, with $u_0, v_0, w_0\geq 0$ on $\bar{\Omega}$, there exists a unique triplet $(u,v,w)$ of nonnegative functions, uniformly bounded in time and classically solving problem \eqref{problem}.
\end{theorem}
\begin{theorem}[The logistic case]\label{MainTheorem1}
	Under the same hypotheses of Theorem \ref{MainTheorem} and $\beta>1$, let $h$ comply with \eqref{h}.
	Then the same claim holds true whenever $\alpha, \gamma >0$, $m_1, m_2, m_3 \in \R$, and $f$ and $g$ fulfill \eqref{f} for each of the following cases:
	\vspace{0.1cm}
	
	\begin{minipage}[t]{0.5\textwidth}
		\begin{enumerate}[label=$A_{\arabic*}$)]
			\setlength{\itemsep}{1.5mm}
			\setcounter{enumi}{16}
			\item \label{A17} $\alpha, \gamma \in \left(0, \frac{1}{n}\right]$, $m_1>\mathcal{A'}$;
			\item \label{A18} $\alpha, \gamma \in \left(\frac{1}{n},1\right)$, $m_1>\mathcal{B'}$;
		\end{enumerate}
	\end{minipage}
	\hspace{-4mm}
	\begin{minipage}[t]{0.5\textwidth}
		\begin{enumerate}[label=$A_{\arabic*}$)]
			\setlength{\itemsep}{1.5mm}
			\setcounter{enumi}{18}
			\item \label{A19} $\alpha \in \left(0,\frac{1}{n}\right], \gamma \in \left(\frac{1}{n},1\right)$, $m_1>\mathcal{C'}$;
			\item \label{A20} $\alpha \in \left(\frac{1}{n},1\right), \gamma \in \left(0,\frac{1}{n}\right]$, $m_1>\mathcal{C'}^t$.
		\end{enumerate}
	\end{minipage}

	\medskip
	\noindent
\end{theorem}
All these results are put together into Table \ref{Table_ResultUnified}. 
\newcolumntype{H}{>{\setbox0=\hbox\bgroup}c<{\egroup}@{}}
\setlength\extrarowheight{4.8pt}
\begin{table}[h!]
	\makegapedcells
	\centering 
	\begin{tabular}{ccccccccccH}
		\hline
		&&&&\multicolumn{3}{c}{{\bf{The non-logistic case}} -- Theorem \ref{MainTheorem}}&&&&\\
		\hline
		$m_2$ &$m_3$&$\alpha$&\multicolumn{1}{c|}{$\gamma$}&\multicolumn{1}{c|}{$m_1$}&$\chi$&$\xi$ 
		\\
		\hline
		$\R$ &$\R$&$(0,\frac{1}{n}]$&\multicolumn{1}{c|}{$(0,\frac{1}{n}]$}&\multicolumn{1}{c|}{$>\mathcal{A}$}&$\R^+$&\multicolumn{1}{c|}{$\R^+$}& 
		\\
		$\R$ &$\R$&$(\frac{1}{n},\frac{2}{n})$&\multicolumn{1}{c|}{$(\frac{1}{n},\frac{2}{n})$}&\multicolumn{1}{c|}{$>\mathcal{B}$}&$\R^+$&\multicolumn{1}{c|}{$\R^+$}
		\\
		$\R$ &$\R$&$[\frac{2}{n},1)$&\multicolumn{1}{c|}{$[\frac{2}{n},1)$}&\multicolumn{1}{c|}{$>\mathcal{C}$}&$\R^+$&\multicolumn{1}{c|}{$\R^+$}& 
\\  
		$\R$ &$\R$&$[\frac{2}{n},1]$&\multicolumn{1}{c|}{$[\frac{2}{n},1)$}&\multicolumn{1}{c|}{$>\mathcal{D}$}&$\R^+$&\multicolumn{1}{c|}{$\R^+$}&
		\\  
		$\R$ &$\R$&$[\frac{2}{n},1)$&\multicolumn{1}{c|}{$[\frac{2}{n},1]$}&\multicolumn{1}{c|}{$>\mathcal{E}$}&$\R^+$&\multicolumn{1}{c|}{$\R^+$}& 
		\\  
		$\R$ &$\R$&$[\frac{2}{n},1]$&\multicolumn{1}{c|}{$[\frac{2}{n},1]$}&\multicolumn{1}{c|}{$>\mathcal{F}$}&$\R^+$&\multicolumn{1}{c|}{$\R^+$}& 
		\\  
		$\R$ &$\R$&$(0,\frac{1}{n}]$&\multicolumn{1}{c|}{$(\frac{1}{n},\frac{2}{n})$}&\multicolumn{1}{c|}{$>\mathcal{G}$}&$\R^+$&\multicolumn{1}{c|}{$\R^+$}&
\\  
		$\R$ &$\R$&$(0,\frac{1}{n}]$&\multicolumn{1}{c|}{$[\frac{2}{n},1)$}&\multicolumn{1}{c|}{$>\mathcal{H}$}&$\R^+$&\multicolumn{1}{c|}{$\R^+$}&
		\\ 
		$\R$ &$\R$&$(0,\frac{1}{n}]$&\multicolumn{1}{c|}{$[\frac{2}{n},1]$}&\multicolumn{1}{c|}{$>\mathcal{I}$}&$\R^+$&\multicolumn{1}{c|}{$\R^+$}&
		\\ 
		$\R$ &$\R$&$(\frac{1}{n},\frac{2}{n})$&\multicolumn{1}{c|}{$(0,\frac{1}{n}]$}&\multicolumn{1}{c|}{$>\mathcal{G}^t$}&$\R^+$&\multicolumn{1}{c|}{$\R^+$}&
\\ 
		$\R$&$\R$&$(\frac{1}{n},\frac{2}{n})$&\multicolumn{1}{c|}{$[\frac{2}{n},1)$}&\multicolumn{1}{c|}{$>\mathcal{J}$}&$\R^+$&\multicolumn{1}{c|}{$\R^+$}&
\\ 
		$\R$ &$\R$&$(\frac{1}{n},\frac{2}{n})$&\multicolumn{1}{c|}{$[\frac{2}{n},1]$}&\multicolumn{1}{c|}{$>\mathcal{K}$}&$\R^+$&\multicolumn{1}{c|}{$\R^+$}&
		\\ 
		$\R$ &$\R$&$[\frac{2}{n},1)$&\multicolumn{1}{c|}{$(0,\frac{1}{n}]$}&\multicolumn{1}{c|}{$>\mathcal{H}^t$}&$\R^+$&\multicolumn{1}{c|}{$\R^+$}&
		\\ 
		$\R$ &$\R$&$[\frac{2}{n},1]$&\multicolumn{1}{c|}{$(0,\frac{1}{n}]$}&\multicolumn{1}{c|}{$>\mathcal{I}^t$}&$\R^+$&\multicolumn{1}{c|}{$\R^+$}&
		\\
		$\R$&$\R$&$[\frac{2}{n},1)$&\multicolumn{1}{c|}{$(\frac{1}{n},\frac{2}{n})$}&\multicolumn{1}{c|}{$>\mathcal{J}^t$}&$\R^+$&\multicolumn{1}{c|}{$\R^+$}&
		\\
		$\R$&$\R$&$[\frac{2}{n},1]$&\multicolumn{1}{c|}{$(\frac{1}{n},\frac{2}{n})$}&\multicolumn{1}{c|}{$>\mathcal{K}^t$}&$\R^+$&\multicolumn{1}{c|}{$\R^+$}&
		\\
		\hline
		&&&&\multicolumn{3}{c}{{\bf{The logistic case}} -- Theorem \ref{MainTheorem1}}&&&&\\
		\hline
		$m_2$ &$m_3$&$\beta$&$\alpha$&\multicolumn{1}{c|}{$\gamma$}&\multicolumn{1}{c|}{$m_1$}&$\chi$&$\xi$&$k$&$\mu$&\\
		\hline
		$\R$ &$\R$&$>1$&$(0,\frac{1}{n}]$&\multicolumn{1}{c|}{$ (0,\frac{1}{n}]$}&\multicolumn{1}{c|}{$>\mathcal{A'}$}&$\R^+$&$ \R^+$&$ \R$&\multicolumn{1}{c|}{$\R^+$}&
		Th.\ \ref{MainTheorem1}\\
		$\R$ &$\R$&$>1$&$(\frac{1}{n},1)$&\multicolumn{1}{c|}{$ (\frac{1}{n},1)$}&\multicolumn{1}{c|}{$>\mathcal{B'}$}&$\R^+$&$ \R^+$&$ \R$&\multicolumn{1}{c|}{$\R^+$}&
		Th.\  \ref{MainTheorem1}\\
		$\R$ &$\R$&$>1$&$(0,\frac{1}{n}]$&\multicolumn{1}{c|}{$ (\frac{1}{n},1)$}&\multicolumn{1}{c|}{$>\mathcal{C'}$}&$\R^+$&$ \R^+$&$ \R$&\multicolumn{1}{c|}{$\R^+$}&
		Th.\  \ref{MainTheorem1}\\
		$\R$ &$\R$&$>1$&$(\frac{1}{n},1)$&\multicolumn{1}{c|}{$ (0,\frac{1}{n}]$}&\multicolumn{1}{c|}{$>\mathcal{C'}^t$}&$\R^+$&$ \R^+$&$ \R$&\multicolumn{1}{c|}{$\R^+$}&
		Th.\  \ref{MainTheorem1}\\
		\hline  
	\end{tabular}
	\caption{With reference to $\S$\ref{NotazioniSection}, this table collects the ranges of the parameters involved in model \eqref{problem} for which boundedness of its solutions is established for any fixed initial distribution $u_0, v_0$ and $w_0$.}
	\label{Table_ResultUnified}
\end{table}

\begin{remark}[On the validity of Theorem \ref{MainTheorem}]\label{CommentiNonLogistico}
	For $m_i$, with $i=1,2,3$, complying with its related assumptions, Theorem \ref{MainTheorem}
	provides a rather complete picture concerning boundedness of solutions to \eqref{problem}.
	Conversely, as far as the linear diffusion and sensitivities ($m_1=m_2=m_3=1$) model is concerned, it still holds except when $\alpha$ and/or $\gamma$ belong to the intervals $\left[\frac{2}{n},1\right]$ and/or $\left[\frac{2}{n},1\right)$, namely in the assumptions \ref{A3}---\ref{A6}, \ref{A8}, \ref{A9} and \ref{A11}---\ref{A16}. Henceforth, what happens in these situations? In low dimensions or under further restrictions on the data, some cases can be saved. Precisely: 
	\begin{enumerate}[label=\Roman*)]
		\setlength{\itemsep}{0.7mm}
		\item If $\alpha=\gamma=1$ (included in \ref{A6}), boundedness of solutions is ensured by requiring $0<\chi<\frac{1}{5n\|v_0\|_{L^{\infty}(\Omega)}}$ and $0<\xi<\frac{1}{5n\|w_0\|_{L^{\infty}(\Omega)}}$, with $n \geq 2$, as seen in 
		\cite[Theorem 1.1]{FrassuGalvanViglialoro};
		\item \label{NonLogistic2} If $\alpha\in\left(0,\frac{1}{n}\right]$, $\gamma \in \left[\frac{2}{n},1\right]$ (\ref{A9}) or $\alpha \in \left[\frac{2}{n},1\right]$, $\gamma \in \left(0,\frac{1}{n}\right]$ (\ref{A14}), 
		boundedness of solutions is achieved for $\alpha \in (0,1]$, $\gamma=1$ (so for $n=1$), for $\alpha \in \left(0,\frac{1}{2}\right]$, $\gamma=1$ and 
		$0<\xi<K_2(n,\lVert  w_0\rVert_{L^\infty(\Omega)})$ (so for $n=2$), or by symmetry, for $\alpha=1$, $\gamma \in (0,1]$ (so for $n=1$), for $\alpha=1$, $\gamma \in \left(0,\frac{1}{2}\right]$ and $0<\chi<K_1(n,\lVert  v_0\rVert_{L^\infty(\Omega)})$ (so for $n=2$), respectively;
		\item \label{NonLogistic3} If $\alpha\in\left(\frac{1}{n},\frac{2}{n}\right)$, $\gamma \in \left[\frac{2}{n},1\right]$ (\ref{A12}) or $\alpha \in \left[\frac{2}{n},1\right]$, 
		$\gamma \in \left(\frac{1}{n},\frac{2}{n}\right)$ (\ref{A16}), Theorem ~\ref{MainTheorem} still holds for $n=1$, $\alpha\in(1,2)$, $\gamma=1$, for $n=2$, $\alpha\in\left(\frac{1}{2},1\right)$, 
		$\gamma=1$ and $0<\xi<\tilde{K}_2(n,\lVert  w_0\rVert_{L^\infty(\Omega)})$, or for $n=1$, $\alpha=1$, $\gamma \in (1,2)$, for $n=2$, $\alpha=1$, $\gamma \in \left(\frac{1}{2},1\right)$ and $0<\chi<\tilde{K}_1(n,\lVert  v_0\rVert_{L^\infty(\Omega)})$, respectively.
	\end{enumerate}
Apparently, the remaining scenarios cannot be managed; some details on this discussion will be given in Remark \ref{DettagliCasoLineareNonLogistico}.
\end{remark}
\begin{remark}[On the validity of Theorem \ref{MainTheorem1}]\label{CommentiLogistico}
	If for the non logistic case Table \ref{Table_ResultUnified} provides an exhaustive scenario, for the logistic one we discuss separately the linear and nonlinear situations.
	
	In particular, for $m_1=m_2=m_3=1$, it is seen that Theorem \ref{MainTheorem1} still applies for $\beta>2$ and $\alpha, \gamma \in (0,1)$. Which is the situation in those cases where the limit values of $\beta$ and/or $\alpha,\gamma$ are considered? Theorem \ref{MainTheorem1} applies (hints are given in  Remark \ref{DettagliCasoLineareLogistico} at the end of the article) whenever
	\begin{enumerate}[label=\roman*)] 
		\setlength{\itemsep}{0.7mm}
		\item \label{logistico1} $\beta>2$ and $\alpha,\gamma\in (0,1]$;
				\item \label{logistico2} $\beta=2, \alpha, \gamma \in (\frac{1}{n},1]$ and
		$\mu > K(n) (\chi^2 \lVert \chi v_0\rVert_{L^\infty(\Omega)}^\frac{4}{n}+ \xi^2 \lVert \xi w_0\rVert_{L^\infty(\Omega)}^\frac{4}{n})$, with $K(n)>0$; 
		\item \label{logistico3} $\beta=2, \alpha \in \left(\frac{1}{n},1\right]$, $\gamma \in \left(0,\frac{1}{n}\right]$ and 
		$\mu > K_1(n) \chi^2 \lVert \chi v_0\rVert_{L^\infty(\Omega)}^\frac{4}{n}$, with  $K_1(n)>0$;
		\item \label{logistico4} $\beta=2,\alpha \in \left(0,\frac{1}{n}\right]$, $\gamma \in \left(\frac{1}{n},1\right]$ and 
		$\mu > K_2(n) \xi^2 \lVert \xi w_0\rVert_{L^\infty(\Omega)}^\frac{4}{n}$, with  $K_2(n)>0$.
	\end{enumerate}
	On the other hand, for the nonlinear diffusion and sensitivities case, it continues valid when either $\beta>2$ and 
	\begin{enumerate}[resume, label=\roman*)]
		\setlength{\itemsep}{0.7mm}
		\item \label{logistico5} $\alpha, \gamma \in \left(\frac{1}{n},1\right]$ and \;  $m_1>\mathcal{B'}$;
		\item \label{logistico6} $\alpha \in \left(0,\frac{1}{n}\right], \gamma \in \left(\frac{1}{n},1\right]$ and \; $m_1>\mathcal{C'}$;
		\item \label{logistico7} $\alpha \in \left(\frac{1}{n},1\right], \gamma \in \left(0,\frac{1}{n}\right]$ and \; $m_1>\mathcal{C'}^t$;
	\end{enumerate}
	or $\beta=2$ and
	\begin{enumerate}[resume, label=\roman*)]
		\setlength{\itemsep}{0.7mm}
		\item \label{logistico8} $\alpha, \gamma \in \left(\frac{1}{n},1\right]$, $m_1>\mathcal{B'}$ and
		$\mu > \tilde{K}$, with some constant $\tilde{K}>0$ depending on $n, m_1, m_2, m_3, \chi, \lVert v_0\rVert_{L^\infty(\Omega)}, \xi, \lVert w_0\rVert_{L^\infty(\Omega)}$;
		\item \label{logistico9} $\alpha \in \left(0,\frac{1}{n}\right], \gamma \in \left(\frac{1}{n},1\right]$, $m_1>\mathcal{C'}$ and $\mu > \tilde{K}_1$, with some constant $\tilde{K}_1>0$ depending on $n, m_1, m_2, m_3, \xi, \lVert w_0\rVert_{L^\infty(\Omega)}$;
		\item \label{logistico10} $\alpha \in \left(\frac{1}{n},1\right], \gamma \in \left(0,\frac{1}{n}\right]$, $m_1>\mathcal{C'}^t$ and 
		$\mu > \tilde{K}_2$, with some constant $\tilde{K}_2>0$
		depending on $n, m_1, m_2, m_3, \chi, \lVert v_0\rVert_{L^\infty(\Omega)}$. 
	\end{enumerate}
\end{remark}
\subsection{Technical strategy and structure of the article}
As it will be formally precised below, the mathematical requirements toward boundedness are connected to some \textit{a priori} estimates of $\int_\Omega u^p=\int_\Omega u^p(x,t)dx$, for some $p>1$ and some $t>0$, with $(u,v,w)$ being any given solution to problem \eqref{problem}.  Whereas in \cite{FrassuGalvanViglialoro}, \textit{inter alia}, dealing with the situation where $m_1=m_2=m_3=1$, this issue is addressed by the employment of a functional of the form $\int_\Omega u^p \varphi$, being $ \varphi$ a suitable function of $(v,w)$, apparently for the nonlinear case under investigation the same attempt does not work. In the specific, the functional herein employed is a natural adjustment of that used in many fully-parabolic Keller--Segel systems with two unknowns, and its expression is, for some $p, q, r>1$ properly large and some $t>0$,
\begin{equation}\label{FunctionalOnLocalSol}
	y(t):=\int_\Omega (u+1)^p + \int_\Omega |\nabla v|^{2q} + \int_\Omega |\nabla w|^{2r}. 
\end{equation}
An evolutive analysis of such a functional leads to a crucial absorption differential inequality in time for the functional itself,  and in turn to the desired uniform-in-time bound for $\int_\Omega u^p.$

The rest of the paper is structured as follows. First, in $\S$\ref{LocalSol}  we give some hints concerning the local existence and uniqueness of a classical solution to model \eqref{problem} and some of its main properties. In this same section we give a \textit{boundedness criterion}, establishing globality and boundedness of local solutions from proper \textit{a priori} $L^p$-boundedness. In turn, in $\S$\ref{EstimatesAndProofSection} we focus on the derivation of these bounds, by means of which we can deduce the claims of Theorems \ref{MainTheorem} and \ref{MainTheorem1}.
\section{Existence of local-in-time solutions and basic properties}\label{LocalSol}
Once $\Omega$, $\chi,\xi$, $m_1,m_2,m_3$ and $f, g, h$ are picked as above, with $(u, v, w)$ we will indicate the classical and nonnegative solution of problem \eqref{problem} defined for all $(x,t) \in \bar{\Omega}\times [0,T_{\rm{max}})$, for some finite $T_{\rm{max}}$, and emanating from the nonnegative initial data $(u_0,v_0, w_0)\in (W^{1,\infty}(\Omega))^3$. In particular, $u$, $v$ and $w$ are such that 
\begin{equation}\label{massConservation}
	\hspace{-0.1cm}
	\int_\Omega u(x, t)dx \leq 
	m_0:=\min\left\{m,\ \left(\frac{k_+}{\mu}\right)^{\frac{1}{\beta-1}}|\Omega|\right\} \textrm{ on }  (0,\TM) \textrm{ and }  m=\int_\Omega u_0, 
\end{equation}
and 
\begin{equation}\label{Cg0}
	0 \leq v\leq \lVert v_0\rVert_{L^\infty(\Omega)} \quad \textrm{and}\quad 
	0 \leq w\leq \lVert w_0\rVert_{L^\infty(\Omega)}\quad \textrm{in}\quad \Omega \times (0,\TM).
\end{equation}
Further, globality and boundedness of $(u,v,w)$ (in  the sense of \eqref{ClassicalAndGlobability}) are ensured whenever (\textit{boundedness criterion}, below)  $u\in L^\infty((0,T_{\rm{max}});L^p(\Omega))$, with $p>1$ arbitrarily large, and uniformly with respect $t\in (0,T_{\rm{max}})$: formally, 
\begin{equation}\label{BoundednessCriterio}
	\begin{array}{c}
		\textrm{If for some } \; C>0 \textrm{ and $p=p(n,m_1,m_2,m_3)>1$ arbitrarily large, } \\ \textrm{we have } \displaystyle \int_\Omega u^p \leq C  \textrm{ on }  (0,\TM), \textrm{ then } (u,v,w) \in (L^\infty((0, \infty);L^{\infty}(\Omega)))^3.
	\end{array}
\end{equation}
We do not prove these basic statements, nor dedicate any lemma; we understand that the details in \cite[$\S$2]{FrassuGalvanViglialoro} and \cite[Appendix A]{TaoWinkParaPara}, which take into considerations also relations \eqref{massConservation} and \eqref{Cg0}, are sufficient in this regard.  Conversely, we spend some words regarding estimate  \eqref{massConservation}. When $h\equiv 0$, it immediately follows by integrating over $\Omega$ the first equation of \eqref{problem} and it is the well-known mass conservation property. In the presence of the logistic terms $h$ as in \eqref{h}, oppositely, also an application of the H\"{o}lder inequality has to be invoked: precisely for $k_+=\max\{k,0\}$ and  
for all $t \in (0,\TM)$
\begin{align*}
	\frac{d}{dt} \int_\Omega u = \int_\Omega h(u) =k \int_\Omega u - \mu \int_\Omega u^{\beta} \leq k_+ \int_\Omega u - \frac{\mu}{|\Omega|^{\beta-1}} \left(\int_\Omega u\right)^{\beta}, 
\end{align*}
so concluding by virtue of an ODI-comparison argument.

Crucial in our computations, beyond the above derivations, are as well some uniform bounds for $\|v(\cdot,t)\|_{W^{1,s}(\Omega)}$ and $\|w(\cdot,t)\|_{W^{1,s}(\Omega)}$, with $s\geq 1$. In this sense, the following lemma is a cornerstone.
\begin{lemma}\label{LocalV}   
For some $c_0, c_1>0$, we have that  $v$ and $w$ comply with
\begin{equation}\label{Cg}
\int_\Omega |\nabla v(\cdot, t)|^s\leq c_0 \quad \textrm{on } \,  (0,\TM)
\begin{cases}
\; \textrm{for all } s \in [1,\infty) & \textrm{if } \alpha \in \left(0, \frac{1}{n}\right],\\
\;  \textrm{for all } s \in \left[1, \frac{n}{(n\alpha-1)}\right) & \textrm{if } \alpha \in \left(\frac{1}{n},1\right],
\end{cases}  
\end{equation}
and
\begin{equation}\label{Cg1}
\int_\Omega |\nabla w(\cdot, t)|^s\leq c_1 \quad \textrm{on } \,  (0,\TM)
\begin{cases}
\; \textrm{for all } s \in [1,\infty) & \textrm{if } \gamma \in \left(0, \frac{1}{n}\right],\\
\;  \textrm{for all } s \in \left[1, \frac{n}{(n\gamma-1)}\right) & \textrm{if } \gamma \in \left(\frac{1}{n},1\right].
\end{cases}  
\end{equation}
\begin{proof}
Fixed $\alpha, \gamma \in (0,1]$, it is possible to find $\rho, \rho_1 >\frac{1}{2}$ such that for all $s \in \left[\frac{1}{\alpha},\frac{n}{(n\alpha-1)_+}\right)$ and  $s \in \left[\frac{1}{\gamma},\frac{n}{(n\gamma-1)_+}\right)$, we have 
$\frac{1}{2}<\rho <1-\frac{n}{2}\big(\alpha-\frac{1}{s}\big)$ and 
$\frac{1}{2}<\rho_1 <1-\frac{n}{2}\big(\gamma-\frac{1}{s}\big)$, respectively. From $1-\rho-\frac{n}{2}\big(\alpha-\frac{1}{s}\big)>0$ and $1-\rho_1-\frac{n}{2}\big(\gamma-\frac{1}{s}\big)>0$, the claims follow invoking properties related to the Neumann heat semigroup; details can be found in \cite[Lemma 5.1]{FrassuCorViglialoro}.
\end{proof}
\end{lemma}
We will also make use of these technical results.
\begin{lemma}\label{LemmaCoefficientAiAndExponents}  
Let $n\in \N$, with $n\geq 2$, $m_1>\frac{n-2}{n}$, $m_2,m_3\in \R$  and $\alpha, \gamma \in (0,1]$. Then there is $s \in [1,\infty)$, such that for proper 
$p,q, r\in[1,\infty)$, $\theta$ and $\theta'$, $\tilde{\theta}$ and $\tilde{\theta}'$, $\mu$ and $\mu'$, $\tilde{\mu}$ and $\tilde{\mu}'$ conjugate exponents, we have that  
\begin{align}
a_1&= \frac{\frac{m_1+p-1}{2}\left(1-\frac{1}{(p+2m_2-m_1-1)\theta}\right)}{\frac{m_1+p-1}{2}+\frac{1}{n}-\frac{1}{2}},            
&  a_2&=\frac{q\left(\frac{1}{s}-\frac{1}{2\theta'}\right)}{\frac{q}{s}+\frac{1}{n}-\frac{1}{2}},  \nonumber \\ 
a_3 &= \frac{\frac{m_1+p-1}{2}\left(1-\frac{1}{2\alpha\mu}\right)}{\frac{m_1+p-1}{2}+\frac{1}{n}-\frac{1}{2}},  
&  a_4  &= \frac{q\left(\frac{1}{s}-\frac{1}{2(q-1)\mu'}\right)}{\frac{q}{s}+\frac{1}{n}-\frac{1}{2}},\nonumber\\ 
\kappa_1 &  =\frac{\frac{p}{2}\left(1- \frac{1}{p}\right)}{\frac{m_1+p-1}{2}+\frac{1}{n}-\frac{1}{2}} \nonumber, & \kappa_2 & =  \frac{q - \frac{1}{2}}{q+\frac{1}{n}-\frac{1}{2}}, 
\end{align}
and
\begin{align}
\tilde{a}_1&= \frac{\frac{m_1+p-1}{2}\left(1-\frac{1}{(p+2m_3-m_1-1)\tilde{\theta}}\right)}{\frac{m_1+p-1}{2}+\frac{1}{n}-\frac{1}{2}},            
&  \tilde{a}_2&=\frac{r\left(\frac{1}{s}-\frac{1}{2\tilde{\theta}'}\right)}{\frac{r}{s}+\frac{1}{n}-\frac{1}{2}},  \nonumber \\ 
\tilde{a}_3 &= \frac{\frac{m_1+p-1}{2}\left(1-\frac{1}{2\gamma\tilde{\mu}}\right)}{\frac{m_1+p-1}{2}+\frac{1}{n}-\frac{1}{2}},  
&  \tilde{a}_4  &= \frac{r\left(\frac{1}{s}-\frac{1}{2(r-1)\tilde{\mu}'}\right)}{\frac{r}{s}+\frac{1}{n}-\frac{1}{2}},\nonumber\\ 
\kappa_3 & =  \frac{r - \frac{1}{2}}{r+\frac{1}{n}-\frac{1}{2}}, \nonumber
\end{align}
belong to the interval $(0,1)$. If, additionally, 
\begin{equation}\label{Restrizionem1-m2-alphaPiccolo}
\alpha \in \left(0,\frac{1}{n}\right] \; \textrm{and}\quad m_1>m_2-\frac{1}{n}, \quad \gamma \in \left(0,\frac{1}{n}\right] \; \textrm{and}\quad m_1>m_3-\frac{1}{n},
\end{equation}
\begin{equation}\label{Restrizionem1-m2-alphaGrande}
\alpha \in \left(\frac{1}{n},\frac{2}{n}\right) \; \textrm{and}\quad m_1>m_2-\frac{2}{n}+\alpha, \quad \gamma \in \left(\frac{1}{n},\frac{2}{n}\right) \; \textrm{and}\quad m_1>m_3-\frac{2}{n}+\gamma, 
\end{equation}
or
\begin{equation}\label{Restrizionem1-m2-alphaGrandeBis}
\alpha \in \left[\frac{2}{n},1\right] \; \textrm{and}\quad m_1>m_2 + \frac{n\alpha-2}{n\alpha-1}, \quad \gamma \in \left[\frac{2}{n},1\right] \; \textrm{and} \quad m_1>m_3 + \frac{n\gamma-2}{n\gamma-1},
\end{equation}
these further relations hold true: 
\begin{align*}
&\beta_1 + \gamma_1 =\frac{p+2m_2-m_1-1}{m_1+p-1}a_1+\frac{1}{q}a_2\in (0,1) \;\textrm{ and }\\
&\beta_2 + \gamma_2= \frac{2 \alpha}{m_1+p-1}a_3+\frac{q-1}{q}a_4 \in (0,1),
\end{align*}
and
\begin{align*}
&\tilde{\beta}_1 + \tilde{\gamma}_1 =\frac{p+2m_3-m_1-1}{m_1+p-1}\tilde{a}_1+\frac{1}{r}\tilde{a}_2\in (0,1) \;\textrm{ and }\\
&\tilde{\beta}_2 + \tilde{\gamma}_2= \frac{2 \gamma }{m_1+p-1}\tilde{a}_3+\frac{r-1}{r}\tilde{a}_4 
\in (0,1).
\end{align*}
Finally, the relations involving the sum of $\beta_1$ and $\gamma_1$, $\beta_2$ and $\gamma_2$, $\tilde{\beta}_1$ and $\tilde{\gamma}_1$, $\tilde{\beta}_2$ and $\tilde{\gamma}_2$ still hold true in each one of the following cases:

	\begin{itemize}
		\setlength{\itemsep}{0.7mm}
		\item [$\triangleright$] $\alpha \in \left(0,\frac{1}{n}\right] \; \textrm{and}\quad m_1>m_2-\frac{1}{n}, \qquad \quad \gamma \in \left(\frac{1}{n}, \frac{2}{n}\right) \; \textrm{and}\quad 
		m_1>m_3-\frac{2}{n}+\gamma$,
		\item [$\triangleright$] $\alpha \in \left(0,\frac{1}{n}\right] \; \textrm{and}\quad m_1>m_2-\frac{1}{n}, \qquad \quad \gamma \in \left[\frac{2}{n},1\right] \; \textrm{and}\quad 
		m_1>m_3 + \frac{n\gamma-2}{n\gamma-1}$,
		\item [$\triangleright$] $\alpha \in \left(\frac{1}{n},\frac{2}{n}\right) \; \textrm{and}\quad m_1>m_2-\frac{2}{n}+\alpha, \quad \gamma \in \left(0,\frac{1}{n}\right] \; \textrm{and}\quad 
		m_1>m_3-\frac{1}{n}$,
		\item [$\triangleright$] $\alpha \in \left(\frac{1}{n},\frac{2}{n}\right) \; \textrm{and}\quad m_1>m_2-\frac{2}{n}+\alpha, \quad \gamma \in \left[\frac{2}{n},1\right] \; \textrm{and}\quad 
		m_1>m_3+ \frac{n\gamma-2}{n\gamma-1}$,
		\item [$\triangleright$] $\alpha \in \left[\frac{2}{n},1\right] \; \textrm{and}\quad m_1>m_2+\frac{n\alpha-2}{n\alpha-1}, \quad \; \; \gamma \in \left(0,\frac{1}{n}\right] \; \textrm{and}\quad 
		m_1>m_3 -\frac{1}{n}$,
		\item [$\triangleright$] $\alpha \in \left[\frac{2}{n},1\right] \; \textrm{and}\quad m_1>m_2+\frac{n\alpha-2}{n\alpha-1}, \quad \; \; \gamma \in \left(\frac{1}{n},\frac{2}{n}\right) \; \textrm{and}\quad 
		m_1>m_3-\frac{2}{n}+\gamma$.
	\end{itemize}
\begin{proof}
For any $s\geq 1$, we put $\theta', \tilde{\theta}'>\max\left\{\frac{n}{2},\frac{s}{2}\right\}$, $\mu>\max\left\{\frac{1}{2\alpha},\frac{n}{2}\right\}$ and 
$\tilde{\mu}>\max\left\{\frac{1}{2\gamma},\frac{n}{2}\right\}$. Thereafter, for 
\begin{equation}\label{Prt_q}
\begin{cases}
q > \max \left\{\frac{n-2}{n}\theta', \frac{s}{2\mu'}+1\right\}, \quad r > \max \left\{\frac{n-2}{n}\tilde{\theta}', \frac{s}{2\tilde{\mu}'}+1\right\} \\
p>\max \left\{2-\frac{2}{n}-m_1,\frac{1}{\theta}-2m_2+m_1+1, \frac{(2m_2-m_1-1)(n-2)\theta-nm_1+n}{n-(n-2)\theta}, \right. \\
\qquad \qquad \; \left. \frac{2\alpha \mu(n-2)}{n} -m_1+1, \frac{1}{\tilde{\theta}}-2m_3+m_1+1,  \right. \\
\qquad \qquad \; \left. 
\frac{(2m_3-m_1-1)(n-2)\tilde{\theta}-nm_1+n}{n-(n-2)\tilde{\theta}}, \frac{2\gamma\tilde{\mu}(n-2)}{n} -m_1+1\right\}, 
\end{cases}
\end{equation} 
it is possible to check that $a_i, \tilde{a}_i, \kappa_2, \kappa_3\in (0,1)$, for any $i=1,2,3,4.$ On the other hand, $\kappa_1\in (0,1)$ also thanks to the assumption 
$m_1>\frac{n-2}{n}.$ 

As to the second part, we consider three cases: $\alpha \in \left(0,\frac{1}{n}\right]$, $\alpha \in \left(\frac{1}{n},\frac{2}{n}\right)$ and 
$\alpha \in \left[\frac{2}{n},1\right].$ 
\begin{itemize}
\item [$\circ$] $\alpha \in \left(0,\frac{1}{n}\right]$. 
For $s>\frac{2\mu'}{2\mu'-1}$ arbitrarily large, consistently with \eqref{Prt_q}, we take $p=q=s$ and $\theta'=s\omega$, for some $\omega>\frac{1}{2}$. Some standard computations entail
\begin{equation*}
0<	\beta_1+\gamma_1=\frac{s+2m_2-m_1-1-\frac{1}{\theta}}{m_1+s-2+\frac{2}{n}}+\frac{2-\frac{1}{\omega}}{s+\frac{2s}{n}},
\end{equation*}
and
\begin{equation*}
0<	\beta_2+\gamma_2=\frac{2\alpha-\frac{1}{\mu}}{m_1+s-2+\frac{2}{n}}+\frac{2s-2-\frac{s}{\mu'}}{s+\frac{2s}{n}}.
\end{equation*}
In light of the above positions, the largeness of $s$ infers $\theta$ arbitrarily close to $1$, in accordance with $\theta'$ large. Further, by choosing $\omega$ approaching $\frac{1}{2}$, continuity arguments imply that  $\beta_1+\gamma_1<1$ whenever restriction \eqref{Restrizionem1-m2-alphaPiccolo} is satisfied, whereas $\beta_2+\gamma_2<1$ comes from $\mu>\frac{n}{2}.$
\item  [$\circ$] $\alpha \in \left(\frac{1}{n},\frac{2}{n}\right).$ First let $s$ be arbitrarily close to $\frac{n}{n\alpha-1}$ and let $q=\frac{p}{2}$ such that \eqref{Prt_q} is accomplished. Then, it holds that $\max\left\{\frac{s}{2},\frac{n}{2}\right\}=\frac{s}{2}$, so that restriction on $\theta'$ (see above) reads $\theta'>\frac{s}{2}$. Subsequently,
\begin{equation*}
0<	\beta_1+\gamma_1=\frac{p+2m_2-m_1-1-\frac{1}{\theta}}{m_1+p-2+\frac{2}{n}}+\frac{2-\frac{s}{\theta'}}{p+\frac{2s}{n}-s},
\end{equation*}
and
\begin{equation*}
0<	\beta_2+\gamma_2=\frac{2\alpha-\frac{1}{\mu}}{m_1+p-2+\frac{2}{n}}+\frac{p-2-\frac{s}{\mu'}}{p+\frac{2s}{n}-s}.
\end{equation*}
Since from $\theta'>\frac{s}{2}$ we have that $\theta'$ approaches $\frac{n}{2(n\alpha-1)}$, an already used reasoning implies that upon enlarging $p$ condition \eqref{Restrizionem1-m2-alphaGrande} yields $\beta_1+\gamma_1<1$. On the other hand,  in order to have $\beta_2+\gamma_2<1$ we have to invoke the above constrain on $\mu$, i.e., $\mu>\frac{1}{2\alpha}.$
\item  [$\circ$] $\alpha \in \left[\frac{2}{n},1\right].$ By considering in the previous case $\theta'>\frac{n}{2}$, we conclude by means of 
\eqref{Restrizionem1-m2-alphaGrandeBis}.
\end{itemize}
By reasoning similarly to what have done before for the range of $\alpha$ and exchanging $\mu'$ with $\tilde{\mu}'$, $q$ with $r$, $\theta'$ with $\tilde{\theta}'$ and $\alpha$ with $\gamma$, we have the claim for the cases $\gamma \in \left(0,\frac{1}{n}\right]$, $\gamma \in \left(\frac{1}{n},\frac{2}{n}\right)$
and $\gamma \in \left[\frac{2}{n},1\right].$ 

The final part is simply obtained by considering permutations of the ranges of $\alpha$ and $\gamma$.
\end{proof}
\end{lemma}
In the concluding part of the paper we will also invoke the next result, by means of which products of powers will be estimated by suitable sums involving their bases and powers of sums controlled by sums of powers. 
\begin{lemma} \label{LemmaIneq}
	Let $a,b,c \geq 0$ and $d_1, d_2>0$ such that $d_1 + d_2 <1$. Then for all $\epsilon >0$ there exists $d>0$ such that
	\[
	a^{d_1}b^{d_2} \leq \epsilon(a+b) +d. 
	\]
	Moreover, for further  $d_3, d_4, d_5>0$, it is possible to find positive $d_6, \hat{d}$ and $\tilde{d}$ such that
	\[
	a^{d_3} + b^{d_4}+c^{d_5} \geq \hat{d}(a+b+c)^{d_6} -\tilde{d}.
	\]
	\begin{proof}
	The proof is based on manipulations of Young’s inequality and some details are available in   
	\cite[Lemma 4.3]{frassuviglialoro1} and \cite[Lemma 3.3]{MarrasViglialoroMathNach}.
	\end{proof}
\end{lemma}

\begin{remark}\label{RemarkOnS}
In view of its importance in the computations, we have to point out that from the above Lemma \ref{LemmaCoefficientAiAndExponents}, the parameter $s$ can be chosen arbitrarily large only when $\alpha, \gamma \in \left(0,\frac{1}{n}\right]$ (this is connected to Lemma  \ref{LocalV}). In particular, as it will be clear later, in such an interval the terms  $\int_{\Omega} (u+1)^{p+2m_2-m_1-1} \vert \nabla v\lvert^2$, 
$\int_{\Omega} (u+1)^{2\alpha} \vert \nabla v\lvert^{2(q-1)}$, $\int_{\Omega} (u+1)^{p+2m_3-m_1-1} \vert \nabla w\lvert^2$ and 
$\int_{\Omega} (u+1)^{2\gamma} \vert \nabla w\lvert^{2(r-1)}$, appearing in our reasoning when dealing with the control of the functional defined in \eqref{FunctionalOnLocalSol}, can be controlled either invoking the Young inequality or the Gagliardo--Nirenberg one. 
\end{remark}
\section{A priori estimates and proof of the theorems}\label{EstimatesAndProofSection}
\subsection{The non-logistic case}\label{NonLog}
In order to exploit the boundedness criterion \eqref{BoundednessCriterio}, let us analyze the behavior of functional defined in \eqref{FunctionalOnLocalSol}, with $p, q, r>1$ properly large. 

In the spirit of Remark \ref{RemarkOnS}, the first steps toward the uniform bound of $\int_\Omega u^p$ will focus on controlling the evolution in time of the functional $y(t)$ by employing the Young inequality.
\begin{lemma}\label{Estim_general_For_u^p_nablav^2qLemma} 
If $m_1,m_2,m_3 \in \R$ comply with $m_1>\max\{2m_2-1, 2m_3-1, \frac{n-2}{n}\}$ or  $m_1>\max\{m_2-\frac{1}{n}, m_3-\frac{1}{n}, \frac{n-2}{n}\}$ or $m_1>\max\{2m_2-1, m_3-\frac{1}{n}, \frac{n-2}{n}\}$ or $m_1>\max\{m_2-\frac{1}{n}, 2m_3-1, \frac{n-2}{n}\}$ whenever 
$\alpha, \gamma \in (0,\frac{1}{n}]$, or $m_1>\max\{2m_2, 2m_3, \frac{n-2}{n}\}$ whenever $\alpha, \gamma \in (\frac{1}{n},1)$, then there exist $p, q, r>1$ such that $(u,v,w)$ satisfies for some $c_{33}, c_{34}, c_{35}, c_{36}>0$ and for all $t \in (0,\TM)$
\begin{align}\label{MainInequality}
\begin{split}
&\frac{d}{dt} \left(\int_\Omega (u+1)^p + \int_\Omega  |\nabla v|^{2q} + \int_\Omega |\nabla w|^{2r} \right) \\
&+ c_{33} \int_\Omega |\nabla |\nabla v|^q|^2 
+ c_{34} \int_\Omega |\nabla |\nabla w|^r|^2 + c_{35} \int_\Omega |\nabla (u+1)^{\frac{m_1+p-1}{2}}|^2 \leq c_{36}.
\end{split}
\end{align}
\begin{proof}
Let $p=q=r>1$ sufficiently large; moreover, in view of Remark \ref{RemarkOnS}, if necessary we are allowed to arbitrarily enlarge these parameters.

For estimates of the term $\frac{d}{dt} \int_\Omega (u+1)^p$, standard testing procedures provide for all $t \in (0,\TM)$
\begin{equation}\label{Estim_1_For_u^p}
\begin{split}
\frac{d}{dt} \int_\Omega (u+1)^p &=\int_\Omega p(u+1)^{p-1}u_t \\
&= -p(p-1) \int_\Omega (u+1)^{p+m_1-3} |\nabla u|^2 \\
&\quad\,+p(p-1)\chi \int_\Omega u(u+1)^{m_2+p-3} \nabla u \cdot \nabla v \\
&\quad\,-p(p-1)\xi \int_\Omega u(u+1)^{m_3+p-3} \nabla u \cdot \nabla w.
\end{split}
\end{equation}
An application of the Young inequality to the second and the third integral in \eqref{Estim_1_For_u^p} give for $\epsilon_1, \epsilon_2>0$ and some positive $c_2,c_3$
\begin{align}\label{Y0}
\begin{split}
&p(p-1)\chi \int_\Omega u (u+1)^{m_2+p-3} \nabla u \cdot \nabla v \\
&\leq \epsilon_1 \int_\Omega (u+1)^{p+m_1-3} |\nabla u|^2 
+ c_2 \int_\Omega (u+1)^{p+2m_2-m_1-1} |\nabla v|^2 \quad \textrm{on }\, (0,\TM),
\end{split}
\end{align}
and  
\begin{align}\label{Y1}  
\begin{split}
&-p(p-1)\xi \int_\Omega u (u+1)^{m_3+p-3} \nabla u \cdot \nabla w \\
&\leq \epsilon_2 \int_\Omega (u+1)^{p+m_1-3} |\nabla u|^2 
+ c_3 \int_\Omega (u+1)^{p+2m_3-m_1-1} |\nabla w|^2 \quad \textrm{on }\, (0,\TM).
\end{split}
\end{align}
\quad \textbf{Case 1}: $\alpha, \gamma \in (0,\frac{1}{n}]$ and $m_1>\max\{2m_2-1, 2m_3-1, \frac{n-2}{n}\}$ or  
$m_1>\max\{m_2-\frac{1}{n}, m_3-\frac{1}{n}, \frac{n-2}{n}\}$ or $m_1>\max\{2m_2-1, m_3-\frac{1}{n}, \frac{n-2}{n}\}$ or 
$m_1>\max\{m_2-\frac{1}{n}, 2m_3-1, \frac{n-2}{n}\}$.
The Young inequality and bound \eqref{Cg} yield for all $t \in (0,\TM)$
\begin{align}\label{Young}  
\begin{split}
c_2 \int_\Omega (u+1)^{p+2m_2-m_1-1} |\nabla v|^2 &\leq  \epsilon_3 \int_\Omega |\nabla v|^{s} + c_4 \int_\Omega (u+1)^{\frac{(p+2m_2-m_1-1)s}{s-2}}\\
&\leq c_4 \int_\Omega (u+1)^{\frac{(p+2m_2-m_1-1)s}{s-2}} + c_5,
\end{split}
\end{align}
with $\epsilon_3 >0$ and some positive $c_4, c_5$. 

Let us now dedicate to the cases $m_1>\max\left\{2m_2-1,\frac{n-2}{n}\right\}$ and $m_1>m_2-\frac{1}{n}$, respectively. 
From $m_1>2 m_2-1$, we have $\frac{(p+2m_2-m_1-1)s}{s-2} < p$, and for every $\epsilon_4>0$, Young's inequality yields some $c_6>0$ entailing 
\begin{equation}\label{Young4}  
c_4 \int_\Omega (u+1)^{\frac{(p+2m_2-m_1-1)s}{s-2}} \leq \epsilon_4 \int_\Omega (u+1)^p + c_6 \quad \textrm{on }\, (0,\TM),
\end{equation}
with $\epsilon_4>0$ and positive $c_6$.
Further, an application of the Gagliardo--Nirenberg inequality and property \eqref{massConservation} yield 
\begin{equation*}\label{Theta_2}
\theta=\frac{\frac{n(m_1+p-1)}{2}\left(1-\frac{1}{p}\right)}{1-\frac{n}{2}+\frac{n(m_1+p-1)}{2}}\in (0,1),
\end{equation*}
so giving for $c_7, c_8>0$
\begin{equation*}
\begin{split}
\int_{\Omega} (u+1)^p&= \|(u+1)^{\frac{m_1+p-1}{2}}\|_{L^{\frac{2p}{m_1+p-1}}(\Omega)}^{\frac{2p}{m_1+p-1}}\\ 
&\leq c_7 \|\nabla (u+1)^{\frac{m_1+p-1}{2}}\|_{L^2(\Omega)}^{\frac{2p}{m_1+p-1}\theta} \|(u+1)^{\frac{m_1+p-1}{2}}\|_{L^{\frac{2}{m_1+p-1}}(\Omega)}^{\frac{2p}{m_1+p-1}(1-\theta)} \\
&\quad\,+ c_7 \|(u+1)^{\frac{m_1+p-1}{2}}\|_{L^{\frac{2}{m_1+p-1}}(\Omega)}^{\frac{2p}{m_1+p-1}}\\
& \leq c_8 \Big(\int_\Omega |\nabla (u+1)^\frac{m_1+p-1}{2}|^2\Big)^{\kappa_1}+ c_8 \quad \textrm{ for all } t \in(0,\TM).
\end{split}
\end{equation*}
Since $\kappa_1 \in (0,1)$ (see Lemma \ref{LemmaCoefficientAiAndExponents}), for any positive $\epsilon_5$  thanks to the Young inequality we arrive for some positive $c_9>0$ at 
 \begin{equation}\label{GN2}
\epsilon_4 \int_\Omega (u+1)^p  \leq \epsilon_5 \int_\Omega |\nabla (u+1)^\frac{m_1+p-1}{2}|^2  + c_9 \quad \textrm{on } (0,\TM).
\end{equation}
Alternatively, we can treat the integral in the left hand side of \eqref{Young4} by applying the Gagliardo--Nirenberg inequality and again bound \eqref{massConservation}, so having
\begin{equation*}  
\begin{split}
&\int_\Omega (u+1)^{\frac{(p+2m_2-m_1-1)s}{s-2}} = \|(u+1)^{\frac{m_1+p-1}{2}}\|_{L^{\frac{2s(p+2m_2-m_1-1)}{(s-2)(m_1+p-1)}}(\Omega)}^{\frac{2s(p+2m_2-m_1-1)}{(s-2)(m_1+p-1)}}\\ 
&\leq c_{10} \|\nabla (u+1)^{\frac{m_1+p-1}{2}}\|_{L^2(\Omega)}^{\frac{2s(p+2m_2-m_1-1)}{(s-2)(m_1+p-1)}\theta_1} \|(u+1)^{\frac{m_1+p-1}{2}}\|_{L^{\frac{2}{m_1+p-1}}(\Omega)}^{\frac{2s(p+2m_2-m_1-1)}{(s-2)(m_1+p-1)}(1-\theta_1)} \\
&\quad\,+ c_{10} \|(u+1)^{\frac{m_1+p-1}{2}}\|_{L^{\frac{2}{m_1+p-1}}(\Omega)}^{\frac{2s(p+2m_2-m_1-1)}{(s-2)(m_1+p-1)}}\\
& \leq c_{11} \Big(\int_\Omega |\nabla (u+1)^\frac{m_1+p-1}{2}|^2\Big)^{\kappa}+ c_{11} \quad \textrm{ for all } t \in(0,\TM),
\end{split}
\end{equation*}
with
\begin{align*}
&\theta_1 = \frac{\frac{m_1+p-1}{2}\left(1-\frac{s-2}{(p+2m_2-m_1-1)s}\right)}{\frac{m_1+p-1}{2}+\frac{1}{n}-\frac{1}{2}} \textrm{ and } 
\kappa = \frac{s(p+2m_2-m_1-1)-(s-2)}{(s-2)(m_1+p-2 + \frac{2}{n})},
\end{align*}
for some $c_{10}, c_{11}>0$. In particular $ \theta_1\in (0,1)$, and from $m_1 > m_2 - \frac{1}{n}$,  since $s$ can be arbitrarily enlarged,  
continuity arguments imply  $\kappa \in (0,1)$; in addition, the Young inequality allow us to rephrase \eqref{Young4} as
\begin{equation}\label{GNY}
c_4 \int_\Omega (u+1)^{\frac{(p+2m_2-m_1-1)s}{s-2}} \leq \epsilon_6 \int_\Omega |\nabla (u+1)^\frac{m_1+p-1}{2}|^2  + c_{12} \quad \textrm{on } (0,\TM),
\end{equation}
with $\epsilon_6 >0$ and some positive $c_{12}$.

Treating in a similar way the second integral on the right-hand side of \eqref{Y1} and exploiting bound \eqref{Cg1} yield 
\begin{align}\label{Young1}  
\begin{split}
&c_3 \int_\Omega (u+1)^{p+2m_3-m_1-1} |\nabla w|^2 \\
&\leq c_{13} \int_\Omega (u+1)^{\frac{(p+2m_3-m_1-1)s}{s-2}} + c_{14} \quad \textrm{ for all } t \in(0,\TM),
\end{split}
\end{align}
with positive $c_{13}, c_{14}$.

Let us now turn our attention to the situation where  $m_1>\max\left\{2m_3-1,\frac{n-2}{n}\right\}$ and $m_1>m_3-\frac{1}{n}$, respectively. In the same flavor as before, from $m_1>2 m_3-1$, we have $\frac{(p+2m_3-m_1-1)s}{s-2} < p$, and for every $\epsilon_7, \epsilon_8>0$, the Young  and the Gagliardo--Nirenberg inequalities yield some $c_{15},c_{16}>0$ entailing 
\begin{align}\label{Yxi}
\begin{split}
c_3 \int_\Omega (u+1)^{\frac{(p+2m_3-m_1-1)s}{s-2}} &\leq \epsilon_7 \int_\Omega (u+1)^p + c_{15} \\
&\leq \epsilon_8 \int_\Omega |\nabla (u+1)^\frac{m_1+p-1}{2}|^2  + c_{16}
\quad \textrm{on } (0,\TM).
\end{split}
\end{align}
On the other hand, by exploiting the condition $m_1>m_3-\frac{1}{n}$ again the Gagliardo--Nirenberg inequality yields 
\begin{equation*}  
\begin{split}
&\int_\Omega (u+1)^{\frac{(p+2m_3-m_1-1)s}{s-2}} = \|(u+1)^{\frac{m_1+p-1}{2}}\|_{L^{\frac{2s(p+2m_3-m_1-1)}{(s-2)(m_1+p-1)}}(\Omega)}^{\frac{2s(p+2m_3-m_1-1)}{(s-2)(m_1+p-1)}}\\ 
&\leq c_{17} \|\nabla (u+1)^{\frac{m_1+p-1}{2}}\|_{L^2(\Omega)}^{\frac{2s(p+2m_3-m_1-1)}{(s-2)(m_1+p-1)}\theta_2} \|(u+1)^{\frac{m_1+p-1}{2}}\|_{L^{\frac{2}{m_1+p-1}}(\Omega)}^{\frac{2s(p+2m_3-m_1-1)}{(s-2)(m_1+p-1)}(1-\theta_2)} \\
&\quad\,+ c_{17} \|(u+1)^{\frac{m_1+p-1}{2}}\|_{L^{\frac{2}{m_1+p-1}}(\Omega)}^{\frac{2s(p+2m_3-m_1-1)}{(s-2)(m_1+p-1)}}\\
& \leq c_{18} \Big(\int_\Omega |\nabla (u+1)^\frac{m_1+p-1}{2}|^2\Big)^{\hat{\kappa}}+ c_{18} \quad \textrm{ for all } t \in(0,\TM),
\end{split}
\end{equation*}
with
\begin{align*}
\theta_2 = \frac{\frac{m_1+p-1}{2}\left(1-\frac{s-2}{(p+2m_3-m_1-1)s}\right)}{\frac{m_1+p-1}{2}+\frac{1}{n}-\frac{1}{2}} \in (0,1)  \textrm{ and } 
\hat{\kappa} = \frac{s(p+2m_3-m_1-1)-(s-2)}{(s-2)(m_1+p-2 + \frac{2}{n})},
\end{align*}
for some $c_{17}, c_{18}>0$, and with $\theta_2, \hat{\kappa} \in (0,1)$. In this way, \eqref{Yxi} becomes
\begin{equation}\label{GNY1}
c_3 \int_\Omega (u+1)^{\frac{(p+2m_3-m_1-1)s}{s-2}} \leq \epsilon_9 \int_\Omega |\nabla (u+1)^\frac{m_1+p-1}{2}|^2  + c_{19} \quad \textrm{on } (0,\TM),
\end{equation}
with $\epsilon_9 >0$ and some positive $c_{19}$.
By plugging estimates \eqref{Y0} and \eqref{Y1} into relation \eqref{Estim_1_For_u^p}, and by relying on bounds \eqref{Young}-\eqref{GN2} and \eqref{Young1}, \eqref{Yxi} (or, alternatively to \eqref{Young4} and \eqref{Yxi}, relations \eqref{GNY} and \eqref{GNY1}), infer for appropriate $\tilde{\epsilon}_1>0$ and proper $c_{20}>0$ for all $t \in (0,\TM)$
\begin{align}\label{ClaimU}
\frac{d}{dt} \int_\Omega (u+1)^p \leq \left(-\frac{4p(p-1)}{(m_1+p-1)^2} + \tilde{\epsilon}_1 \right) \int_\Omega |\nabla (u+1)^{\frac{m_1+p-1}{2}}|^2 
+ c_{20},
\end{align}
where we also have taken into consideration
\begin{equation}\label{GradU}
\int_\Omega (u+1)^{p+m_1-3} |\nabla u|^2 = \frac{4}{(m_1+p-1)^2} \int_\Omega |\nabla (u+1)^{\frac{m_1+p-1}{2}}|^2 \quad  \textrm{on } (0,\TM).
\end{equation}
Now, we treat the terms $\frac{d}{dt} \int_\Omega  |\nabla v|^{2q}$ and $\frac{d}{dt} \int_\Omega  |\nabla w|^{2r}$ of the functional $y(t)$ under the assumption that 
$m_1>\frac{n-2}{n}$.
As to the term $\frac{d}{dt} \int_\Omega  |\nabla v|^{2q}$, reasoning similarly as in \cite[Lemma 5.3]{FrassuCorViglialoro}, we obtain for some 
$c_{21}, c_{22}>0$ for all $t \in (0,\TM)$
\begin{equation}\label{Estim_gradV}
\frac{d}{dt}\int_\Omega  |\nabla v|^{2q}+ q \int_\Omega |\nabla v|^{2q-2} |D^2v|^2 \leq c_{21} \int_\Omega u^{2\alpha} |\nabla v|^{2q-2} +c_{22}. 
\end{equation}
Moreover, Young's inequality and bound \eqref{Cg} give for every arbitrary $\epsilon_{10}, \epsilon_{11}, \epsilon_{12}>0$ and some $c_{23}, c_{24}, c_{25},c_{26}>0$
\begin{equation}\label{Estimat_nablav^2q+2}
\begin{split}
&c_{21} \int_\Omega u^{2\alpha} |\nabla v|^{2q-2} \leq \epsilon_{10} \int_\Omega u^p + c_{23} \int_\Omega |\nabla v|^{\frac{2(q-1)p}{p-2\alpha}}\\
& \leq \epsilon_{10} \int_\Omega (u+1)^p + \epsilon_{11} \int_\Omega |\nabla v|^s + c_{24} \leq \epsilon_{10} \int_\Omega (u+1)^p + c_{25}\\
& \leq \epsilon_{12} \int_\Omega |\nabla (u+1)^{\frac{m_1+p-1}{2}}|^2 + c_{26} \quad \textrm{on } (0,\TM).
\end{split}
\end{equation}
As to the term $\frac{d}{dt} \int_\Omega  |\nabla w|^{2r}$ of the functional $y(t)$,  with bound \eqref{Cg1} in our hands, through similar aforedescribed computations we obtain for all $t \in (0,\TM)$ 
\begin{equation}\label{Estimat_nablaw^2r+2}
\begin{split}
&\frac{d}{dt}\int_\Omega  |\nabla w|^{2r}+ r \int_\Omega |\nabla w|^{2r-2} |D^2 w|^2 \leq c_{27} \int_\Omega u^{2\gamma} |\nabla w|^{2r-2} +c_{28}\\ 
&\leq \epsilon_{13} \int_\Omega u^p + c_{29} \int_\Omega |\nabla w|^{\frac{2(r-1)p}{p-2\gamma}} 
\leq \epsilon_{13} \int_\Omega (u+1)^p + \epsilon_{14} \int_\Omega |\nabla w|^s + c_{30} \\
&\leq \epsilon_{13} \int_\Omega (u+1)^p + c_{31} \leq \epsilon_{15} \int_\Omega |\nabla (u+1)^{\frac{m_1+p-1}{2}}|^2 + c_{32}, 
\end{split}
\end{equation}
with $\epsilon_{13}, \epsilon_{14}, \epsilon_{15}>0$ and some  $c_{27}, c_{28}, c_{29},c_{30}, c_{31}, c_{32}>0$.
Therefore, by inserting relation \eqref{Estimat_nablav^2q+2} into \eqref{Estim_gradV} and adding \eqref{ClaimU} and \eqref{Estimat_nablaw^2r+2}, we have the claim for a proper choice of  $\tilde{\epsilon}_1$ and some positive $c_{33}, c_{34}, c_{35}, c_{36}$, once relations (see \cite[page 17]{FrassuCorViglialoro})
\begin{align}\label{GradV}
\begin{split}
\vert \nabla \lvert \nabla v\rvert^q\rvert^2&=\frac{q^2}{4}\lvert \nabla v \rvert^{2q-4}\vert \nabla \lvert \nabla v\rvert^2\rvert^2=q^2\lvert \nabla v \rvert^{2q-4}\lvert D^2v \nabla v \rvert^2\\
&\leq q^2|\nabla v|^{2q-2} |D^2v|^2\quad \textrm{on }\, (0,\TM), 
\end{split}
\end{align}
and 
\begin{align}\label{GradW}
\begin{split}
\vert \nabla \lvert \nabla w\rvert^r\rvert^2&=\frac{r^2}{4}\lvert \nabla w \rvert^{2r-4}\vert \nabla \lvert \nabla w\rvert^2\rvert^2=r^2\lvert \nabla w \rvert^{2r-4}\lvert D^2w \nabla w \rvert^2\\
&\leq r^2|\nabla w|^{2r-2} |D^2 w|^2 \quad \textrm{on }\, (0,\TM), 
\end{split}
\end{align}
are considered too. 

\quad \textbf{Case 2}: $\alpha, \gamma \in \left(\frac{1}{n},1\right)$ and $m_1>\max\left\{2m_2, 2m_3, \frac{n-2}{n}\right\}$.
Accordingly to Remark \ref{RemarkOnS}, since in this case $s$ has a finite upper bound, a different approach to deal with relations \eqref{Young}, \eqref{Young1}, \eqref{Estimat_nablav^2q+2} and \eqref{Estimat_nablaw^2r+2} has to be used.  In particular, for $\bar{\epsilon}_1>0$ and some $\bar{c}_1>0$ we can estimate relation \eqref{Young} as follows:
\begin{align*}
&c_2 \int_\Omega (u+1)^{p+2m_2-m_1-1} |\nabla v|^2 \\
&\leq \bar{\epsilon}_1 \int_\Omega |\nabla v|^{2(p+1)} + \bar{c}_1 \int_\Omega (u+1)^{\frac{(p+2m_2-m_1-1)(p+1)}{p}} 
\quad \textrm{on }\, (0,\TM).
\end{align*}
Now, if $m_1>2m_2$, then some $p$ sufficiently large infers $\frac{(p+2m_2-m_1-1)(p+1)}{p}< p$, so that for any positive $\bar{\epsilon}_2, \bar{\epsilon}_3$ 
and some $\bar{c}_2, \bar{c}_3>0$ we have
\begin{align*}
&\bar{c}_1 \int_\Omega (u+1)^{\frac{(p+2m_2-m_1-1)(p+1)}{p}} \leq \bar{\epsilon}_2 \int_\Omega (u+1)^p + \bar{c}_2\\
&\leq \bar{\epsilon}_3 \int_\Omega |\nabla (u+1)^{\frac{m_1+p-1}{2}}|^2 + \bar{c}_3 \quad \textrm{for all } t \in (0,\TM),
\end{align*}
where in the last implication we used $m_1>\frac{n-2}{n}$ (as in the previous case).
In a similar way, for $m_1>2m_3$ and $m_1>\frac{n-2}{n}$ we have for any positive $\bar{\epsilon}_4, \bar{\epsilon}_5, \bar{\epsilon}_6$ and some $\bar{c}_4, \bar{c}_5, \bar{c}_6>0$
\begin{equation*}
\begin{split}
&c_3 \int_\Omega (u+1)^{p+2m_3-m_1-1} |\nabla w|^2 \\
&\leq \bar{\epsilon}_4 \int_\Omega |\nabla w|^{2(p+1)} + \bar{c}_4 \int_\Omega (u+1)^{\frac{(p+2m_3-m_1-1)(p+1)}{p}}\\ 
&\leq \bar{\epsilon}_4 \int_\Omega |\nabla w|^{2(p+1)} + \bar{\epsilon}_5 \int_\Omega (u+1)^p + \bar{c}_5\\ 
&\leq  \bar{\epsilon}_4 \int_\Omega |\nabla w|^{2(p+1)} + \bar{\epsilon}_6 \int_\Omega |\nabla (u+1)^{\frac{m_1+p-1}{2}}|^2 + \bar{c}_6 \quad \textrm{on } (0,\TM).
\end{split}
\end{equation*}
Let us focus on the integrals $\int_\Omega |\nabla v|^{\frac{2p(p-1)}{p-2\alpha}}$ and $\int_\Omega |\nabla w|^{\frac{2p(p-1)}{p-2\gamma}}$.
Since $\alpha,\gamma <1$, this implies that $\frac{2p(p-1)}{p-2\alpha}<2(p+1)$ and $\frac{2p(p-1)}{p-2\gamma}<2(p+1)$, and subsequently an application of the Young inequality 
leads to
\begin{equation*}
c_{23} \int_\Omega |\nabla v|^{\frac{2p(p-1)}{p-2\alpha}} \leq \bar{\epsilon}_7 \int_\Omega |\nabla v|^{2(p+1)} + \bar{c}_7 \quad \textrm{for all } t \in (0,\TM),
\end{equation*}
and 
\begin{equation*}
c_{29} \int_\Omega |\nabla w|^{\frac{2p(p-1)}{p-2\gamma}} \leq \bar{\epsilon}_8 \int_\Omega |\nabla w|^{2(p+1)} + \bar{c}_8 \quad \textrm{on } (0,\TM),
\end{equation*}
with $\bar{\epsilon}_7, \bar{\epsilon}_8>0$ and some positive $\bar{c}_7, \bar{c}_8$.
By taking into account \cite[Lemma 2.2]{LankeitWangConsumptLogistic} and bounds \eqref{Cg0}, we get
\begin{equation*}
\int_\Omega |\nabla v|^{2(p+1)} \leq 2 (4p^2+n) \|v_0\|_{L^{\infty}}^2 \int_\Omega |\nabla v|^{2p-2} |D^2 v|^2 \quad \textrm{ for all } t \in (0,\TM),
\end{equation*}
and 
\begin{equation*}
\int_\Omega |\nabla w|^{2(p+1)} \leq 2 (4p^2+n) \|w_0\|_{L^{\infty}}^2 \int_\Omega |\nabla w|^{2p-2} |D^2 w|^2 \quad \textrm{on } (0,\TM);
\end{equation*}
the rest of the proof is an evident adaptation of previous reasoning.

We conclude by observing that this lemma holds in each of the following cases:
\begin{itemize}
	\item [$\rhd$] $\alpha \in \left(0,\frac{1}{n}\right]$ and $\gamma \in \left(\frac{1}{n},1\right)$, whenever $m_1>\max\left\{2m_2-1, 2m_3, \frac{n-2}{n}\right\}$
	or $m_1>\max\left\{m_2-\frac{1}{n}, 2m_3, \frac{n-2}{n}\right\}$
\item [$\rhd$] $\alpha \in \left(\frac{1}{n},1\right)$ and $\gamma \in \left(0,\frac{1}{n}\right]$, whenever $m_1>\max\left\{2m_2, 2m_3-1,\frac{n-2}{n}\right\}$ 
	or $m_1>\max\left\{2m_2, m_3-\frac{1}{n},\frac{n-2}{n}\right\}$.
\end{itemize}
\end{proof}
\end{lemma}
Let us now turn our attention when, as mentioned before, the Gagliardo--Nirenberg inequality is employed. In this case, we can derive information also for 
$\alpha, \gamma= 1$; this is the reason why in Theorem \ref{MainTheorem} we distinguish the situations where the value 1 belongs or not to the interval in question. 
\begin{lemma}\label{Met_GN}
If $m_1,m_2, m_3 \in \R$ and $\alpha, \gamma>0$ are taken accordingly to \eqref{Restrizionem1-m2-alphaPiccolo}, \eqref{Restrizionem1-m2-alphaGrande}, 
\eqref{Restrizionem1-m2-alphaGrandeBis}, then there exist $p, q, r>1$ such that $(u,v,w)$ satisfies a similar inequality as in \eqref{MainInequality}.
\begin{proof}
For $s$, $p$, $q$ and $r$ taken accordingly to Lemma \ref{LemmaCoefficientAiAndExponents} (in particular, $p=q=r$ for $\alpha, \gamma \in \left(0,\frac{1}{n}\right]$, and $q=r=\frac{p}{2}$ for $\alpha, \gamma \in \left(\frac{1}{n},1\right]$), let $\theta, \theta', \tilde{\theta}, \tilde{\theta}', \mu, \mu', \tilde{\mu}, \tilde{\mu}'$, $a_1,a_2, a_3, a_4, \tilde{a}_1, \tilde{a}_2, \tilde{a}_3, \tilde{a}_4$ and $\beta_1, \beta_2, \tilde{\beta}_1, \tilde{\beta}_2, \gamma_1, \tilde{\gamma}_1, \gamma_2, \tilde{\gamma}_2$ be therein defined.

With a view to Lemma \ref{Estim_general_For_u^p_nablav^2qLemma}, by manipulating relation \eqref{Estim_1_For_u^p} and focusing on the inequalities \eqref{Y0}, \eqref{Y1}, \eqref{Estim_gradV} and on the first inequality in \eqref{Estimat_nablaw^2r+2}, a proper $\tilde{\epsilon}_1$ leads to
\begin{equation}\label{Somma}
\begin{split}
&\quad\,\frac{d}{dt} \left(\int_\Omega (u+1)^p + \int_\Omega  |\nabla v|^{2q} + \int_\Omega  |\nabla w|^{2r}\right)\\
&\quad\, + q \int_\Omega |\nabla v|^{2q-2} |D^2v|^2 
+ r \int_\Omega |\nabla w|^{2r-2} |D^2 w|^2\\
&\leq \left(-\frac{4p(p-1)}{(m_1+p-1)^2} + \tilde{\epsilon}_1 \right) \int_\Omega |\nabla (u+1)^{\frac{m_1+p-1}{2}}|^2\\
&\quad\,+  c_2 \int_\Omega (u+1)^{p+2m_2-m_1-1} |\nabla v|^2 + c_{21} \int_\Omega u^{2\alpha} |\nabla v|^{2q-2}\\ 
&\quad\,+  c_3 \int_\Omega (u+1)^{p+2m_3-m_1-1} |\nabla w|^2 + c_{27} \int_\Omega u^{2\gamma} |\nabla w|^{2r-2} \\
&\quad\,+c_{37} \quad \textrm{for all } t \in (0,\TM),
\end{split}
\end{equation} 
for some $c_{37}>0$ (we also used relation \eqref{GradU}).
In this way, we can estimate the last four integrals on the right-hand side of \eqref{Somma} by applying the H\"{o}lder inequality so to have for all $t  \in (0,\TM)$
\begin{align}
&\!\int_{\Omega} (u+1)^{p+2m_2-m_1-1} |\nabla v|^2\leq  \left(\int_{\Omega} (u+1)^{(p+2m_2-m_1-1)\theta}\right)^{\frac{1}{\theta}} 
\left(\int_{\Omega} |\nabla v|^{2 \theta'}\right)^{\frac{1}{\theta'}}, \label{H1} \\
&\!\int_{\Omega} (u+1)^{2\alpha} |\nabla v|^{2q-2}\leq 
\left(\int_{\Omega} (u+1)^{2\alpha\mu}\right)^{\frac{1}{\mu}} \left(\int_{\Omega} |\nabla v|^{2(q-1)\mu'}\right)^{\frac{1}{\mu'}}, \label{H2} 
\end{align}
and
\begin{align} 
&\!\!\!\!\int_{\Omega} (u+1)^{p+2m_3-m_1-1} |\nabla w|^2 
\leq  \left(\int_{\Omega} (u+1)^{(p+2m_3-m_1-1)\tilde{\theta}}\right)^{\frac{1}{\tilde{\theta}}} 
\left(\int_{\Omega} |\nabla w|^{2 \tilde{\theta}'}\right)^{\frac{1}{\tilde{\theta}'}}, \label{H3} \\
&\!\!\!\!\int_{\Omega} (u+1)^{2\gamma} |\nabla w|^{2r-2} \leq \left(\int_{\Omega} (u+1)^{2\gamma\tilde{\mu}}\right)^{\frac{1}{\tilde{\mu}}} \left(\int_{\Omega} |\nabla w|^{2(r-1)\tilde{\mu}'}\right)^{\frac{1}{\tilde{\mu}'}}.  \label{H4}
\end{align}
By invoking the Gagliardo--Nirenberg inequality and bound \eqref{massConservation}, we obtain for some $c_{38}, c_{39}>0$
\begin{align}\label{a1}
\begin{split}
&\left(\int_{\Omega} (u+1)^{(p+2m_2-m_1-1)\theta}\right)^{\frac{1}{\theta}}= \|(u+1)^{\frac{m_1+p-1}{2}}\|_{L^{\frac{2(p+2m_2-m_1-1)}{m_1+p-1}\theta}(\Omega)}^{\frac{2(p+2m_2-m_1-1)}{m_1+p-1}}\\ 
& \leq c_{38} \|\nabla(u+1)^{\frac{m_1+p-1}{2}}\|_{L^2(\Omega)}^{\frac{2(p+2m_2-m_1-1)}{m_1+p-1} a_1} \|(u+1)^{\frac{m_1+p-1}{2}}\|_{L^{\frac{2}{m_1+p-1}}(\Omega)}^{\frac{2(p+2m_2-m_1-1)}{m_1+p-1} (1-a_1)} \\
&\quad\,+ c_{38} \|(u+1)^{\frac{m_1+p-1}{2}}\|_{L^{\frac{2}{m_1+p-1}}(\Omega)}^{\frac{2(p+2m_2-m_1-1)}{m_1+p-1}} \\ 
&\leq c_{39} \left(\int_{\Omega} |\nabla (u+1)^{\frac{m_1+p-1}{2}}|^2\right)^{\beta_1}+ c_{39} \quad \textrm{on } \, (0,\TM),
\end{split}
\end{align}
and for some $c_{40}, c_{41}>0$
\begin{align}\label{tildea1}
\begin{split}
&\left(\int_{\Omega} (u+1)^{(p+2m_3-m_1-1)\tilde{\theta}}\right)^{\frac{1}{\tilde{\theta}}}= \|(u+1)^{\frac{m_1+p-1}{2}}\|_{L^{\frac{2(p+2m_3-m_1-1)}{m_1+p-1}\tilde{\theta}}(\Omega)}^{\frac{2(p+2m_3-m_1-1)}{m_1+p-1}}\\
& \leq c_{40} \|\nabla(u+1)^{\frac{m_1+p-1}{2}}\|_{L^2(\Omega)}^{\frac{2(p+2m_3-m_1-1)}{m_1+p-1} \tilde{a}_1} \|(u+1)^{\frac{m_1+p-1}{2}}\|_{L^{\frac{2}{m_1+p-1}}(\Omega)}^{\frac{2(p+2m_3-m_1-1)}{m_1+p-1} (1-\tilde{a}_1)} \\
&\quad\,+ c_{40} \|(u+1)^{\frac{m_1+p-1}{2}}\|_{L^{\frac{2}{m_1+p-1}}(\Omega)}^{\frac{2(p+2m_3-m_1-1)}{m_1+p-1}} \\ 
&\leq c_{41} \left(\int_{\Omega} |\nabla (u+1)^{\frac{m_1+p-1}{2}}|^2\right)^{\tilde{\beta}_1}+ c_{41} \quad \textrm{ for all } \,t\in(0,\TM).
\end{split}
\end{align}
Moreover, we get for some $c_{42}, c_{43}>0$
\begin{align}\label{a3}
\begin{split}
&\left(\int_{\Omega} (u+1)^{2\alpha\mu}\right)^{\frac{1}{\mu}}= \|(u+1)^{\frac{m_1+p-1}{2}}\|_{L^{\frac{4\alpha\mu}{m_1+p-1}}(\Omega)}^{\frac{4\alpha}{m_1+p-1}}\\ 
& \leq c_{42} \|\nabla(u+1)^{\frac{m_1+p-1}{2}}\|_{L^2(\Omega)}^{\frac{4\alpha}{m_1+p-1}  a_3} \|(u+1)^{\frac{m_1+p-1}{2}}\|_{L^{\frac{2}{m_1+p-1}}(\Omega)}^{\frac{4\alpha}{m_1+p-1} (1-a_3)} \\
&\quad\,+ c_{42} \|(u+1)^{\frac{m_1+p-1}{2}}\|_{L^{\frac{2}{m_1+p-1}}(\Omega)}^{\frac{4\alpha}{m_1+p-1}}\\
&\leq c_{43} \left(\int_{\Omega} |\nabla (u+1)^{\frac{m_1+p-1}{2}}|^2\right)^{\beta_2}+ c_{43}
\quad \textrm{on }\, (0,\TM),
\end{split}
\end{align}
and for some $c_{44}, c_{45}>0$
\begin{align}\label{tildea3}
\begin{split}
&\left(\int_{\Omega} (u+1)^{2\gamma\tilde{\mu}}\right)^{\frac{1}{\tilde{\mu}}}= \|(u+1)^{\frac{m_1+p-1}{2}}\|_{L^{\frac{4\gamma\tilde{\mu}}{m_1+p-1}}(\Omega)}^{\frac{4\gamma}{m_1+p-1}}\\ 
& \leq c_{44} \|\nabla(u+1)^{\frac{m_1+p-1}{2}}\|_{L^2(\Omega)}^{\frac{4\gamma}{m_1+p-1}  \tilde{a}_3} \|(u+1)^{\frac{m_1+p-1}{2}}\|_{L^{\frac{2}{m_1+p-1}}(\Omega)}^{\frac{4\gamma}{m_1+p-1} (1-\tilde{a}_3)} \\
&\quad\,+ c_{44} \|(u+1)^{\frac{m_1+p-1}{2}}\|_{L^{\frac{2}{m_1+p-1}}(\Omega)}^{\frac{4\gamma}{m_1+p-1}}\\
&\leq c_{45} \left(\int_{\Omega} |\nabla (u+1)^{\frac{m_1+p-1}{2}}|^2\right)^{\tilde{\beta}_2}+ c_{45}
\quad \textrm{for all }\, t \in (0,\TM).
\end{split}
\end{align}
In a similar way, we can again apply the Gagliardo--Nirenberg inequality and bound \eqref{Cg} and get for some $c_{46}, c_{47}>0$
\begin{align}\label{a2}
\begin{split}
&\left(\int_{\Omega} |\nabla v|^{2 \theta'}\right)^{\frac{1}{\theta'}}=\| |\nabla v|^q \|_{L^{\frac{2 \theta'}{q}}(\Omega)}^{\frac{2}{q}}\\
&\leq c_{46} \|\nabla |\nabla v|^q\|_{L^2(\Omega)}^{\frac{2}{q}a_2} \| |\nabla v|^q\|_{L^{\frac{s}{q}}(\Omega)}^{\frac{2}{q}(1-a_2)} + c_{46} \| |\nabla v|^q\|_{L^{\frac{s}{q}}(\Omega)}^{\frac{2}{q}}\\ 
&\leq c_{47} \left(\int_{\Omega} |\nabla |\nabla v|^q|^2 \right)^{\gamma_1}+ c_{47}
\quad \textrm{on } (0, \TM),
\end{split}
\end{align}
and for some $c_{48}, c_{49}>0$
\begin{equation}\label{a4}
\begin{split}
&\left(\int_{\Omega} |\nabla v|^{2(q-1) \mu'}\right)^{\frac{1}{\mu'}}=\| |\nabla v|^q \|_{L^{\frac{2(q-1)}{q}\mu'}(\Omega)}^{\frac{2(q-1)}{q}} 
\\
&\leq c_{48} \|\nabla |\nabla v|^q\|_{L^2(\Omega)}^{\frac{2(q-1)}{q}  a_4} \| |\nabla v|^q\|_{L^{\frac{s}{q}}(\Omega)}^{\frac{2(q-1)}{q}  (1-a_4)}+ c_{48} \| |\nabla v|^q\|_{L^{\frac{s}{q}}(\Omega)}^{\frac{2(q-1)}{q}} \\
&\leq c_{49} \left(\int_{\Omega} |\nabla |\nabla v|^q|^2 \right)^{\gamma_2}+ c_{49} \quad \textrm{for all } t\in (0,\TM).
\end{split}
\end{equation}
Finally, an application of the Gagliardo--Nirenberg inequality and bound \eqref{Cg1} imply for some $c_{50}, c_{51}>0$
\begin{equation}\label{tildea2}
\begin{split}
&\left(\int_{\Omega} |\nabla w|^{2 \tilde{\theta}'}\right)^{\frac{1}{\tilde{\theta}'}}=\| |\nabla w|^r \|_{L^{\frac{2 \tilde{\theta}'}{r}}(\Omega)}^{\frac{2}{r}}\\
&\leq c_{50} \|\nabla |\nabla w|^r\|_{L^2(\Omega)}^{\frac{2}{r} \tilde{a}_2} \| |\nabla w|^r\|_{L^{\frac{s}{r}}(\Omega)}^{\frac{2}{r}(1-\tilde{a}_2)}+ c_{50} \| |\nabla w|^r\|_{L^{\frac{s}{r}}(\Omega)}^{\frac{2}{r}} \\
& \leq c_{51} \left(\int_{\Omega} |\nabla |\nabla w|^r|^2 \right)^{\tilde{\gamma}_1}+ c_{51} \quad \textrm{on } (0, \TM),
\end{split}
\end{equation}
and for some $c_{52}, c_{53}>0$
\begin{equation}\label{tildea4}
\begin{split}
&\left(\int_{\Omega} |\nabla w|^{2(r-1) \tilde{\mu}'}\right)^{\frac{1}{\tilde{\mu}'}}=\| |\nabla w|^r \|_{L^{\frac{2(r-1)}{r}\tilde{\mu}'}(\Omega)}^{\frac{2(r-1)}{r}} 
\\
&\leq c_{52} \|\nabla |\nabla w|^r\|_{L^2(\Omega)}^{\frac{2(r-1)}{r}  \tilde{a}_4} \| |\nabla w|^r\|_{L^{\frac{s}{r}}(\Omega)}^{\frac{2(r-1)}{r}  (1-\tilde{a}_4)}+ c_{52} \| |\nabla w|^r\|_{L^{\frac{s}{r}}(\Omega)}^{\frac{2(r-1)}{r}} \\
&\leq c_{53} \left(\int_{\Omega} |\nabla |\nabla w|^r|^2 \right)^{\gamma_2}+ c_{53} \quad \textrm{for all } t\in (0,\TM).
\end{split}
\end{equation}
By plugging \eqref{H1}, \eqref{H2}, \eqref{H3} and \eqref{H4} into \eqref{Somma} and taking into account \eqref{a1}, \eqref{tildea1}, \eqref{a2}, \eqref{tildea2},
\eqref{a3}, \eqref{tildea3}, \eqref{a4}, \eqref{tildea4}, once inequalities \eqref{GradV} and \eqref{GradW} are considered we can derive for a positive suitable $\tilde{\epsilon}_1$ the following estimate, valid for certain $c_{54}, c_{55}, c_{56}, c_{57}, c_{58}, c_{59}, c_{60}, c_{61}>0$ and for all $(0,\TM)$:
\begin{equation}\label{Somma1}
\begin{split}
&\quad\,\frac{d}{dt} \left(\int_\Omega (u+1)^p + \int_\Omega  |\nabla v|^{2q} + \int_\Omega  |\nabla w|^{2r}\right) \\
&\quad\,+ c_{54} \int_\Omega |\nabla |\nabla v|^q|^2 
+ c_{55} \int_\Omega |\nabla |\nabla w|^r|^2  + c_{56} \int_\Omega |\nabla (u+1)^{\frac{m_1+p-1}{2}}|^2 \\
&\leq c_{57} \left(\int_{\Omega} |\nabla (u+1)^{\frac{m_1+p-1}{2}}|^2\right)^{\beta_1}  
\left(\int_{\Omega} |\nabla |\nabla v|^q|^2 \right)^{\gamma_1} \\
&\quad\,+ c_{57} \left(\int_{\Omega} |\nabla (u+1)^{\frac{m_1+p-1}{2}}|^2\right)^{\beta_1}+ c_{57} \left(\int_{\Omega} |\nabla |\nabla v|^q|^2 \right)^{\gamma_1} \\
&\quad\,+ c_{58} \left(\int_{\Omega} |\nabla (u+1)^{\frac{m_1+p-1}{2}}|^2\right)^{\beta_2} \left(\int_{\Omega} |\nabla |\nabla v|^q|^2 \right)^{\gamma_2}\\ 
&\quad\,+ c_{58} \left(\int_{\Omega} |\nabla (u+1)^{\frac{m_1+p-1}{2}}|^2\right)^{\beta_2} + c_{58} \left(\int_{\Omega} |\nabla |\nabla v|^q|^2 \right)^{\gamma_2}\\
&\quad\,+c_{59} \left(\int_{\Omega} |\nabla (u+1)^{\frac{m_1+p-1}{2}}|^2\right)^{\tilde{\beta}_1}  
\left(\int_{\Omega} |\nabla |\nabla w|^r|^2 \right)^{\tilde{\gamma}_1} \\
&\quad\,+ c_{59} \left(\int_{\Omega} |\nabla (u+1)^{\frac{m_1+p-1}{2}}|^2\right)^{\tilde{\beta}_1}+ c_{59} \left(\int_{\Omega} |\nabla |\nabla w|^r|^2 \right)^{\tilde{\gamma}_1}\\ 
&\quad\,+ c_{60} \left(\int_{\Omega} |\nabla (u+1)^{\frac{m_1+p-1}{2}}|^2\right)^{\tilde{\beta}_2} \left(\int_{\Omega} |\nabla |\nabla w|^r|^2 \right)^{\tilde{\gamma}_2}\\ 
&\quad\,+ c_{60} \left(\int_{\Omega} |\nabla (u+1)^{\frac{m_1+p-1}{2}}|^2\right)^{\tilde{\beta}_2} + c_{60} \left(\int_{\Omega} |\nabla |\nabla w|^r|^2 \right)^{\tilde{\gamma}_2} 
+ c_{61}.
\end{split}
\end{equation} 
Since by Lemma \ref{LemmaCoefficientAiAndExponents} we have that $\beta_1 + \gamma_1 <1$, $\beta_2 + \gamma_2 <1$, $\tilde{\beta}_1 + \tilde{\gamma}_1 <1$ and 
$\tilde{\beta}_2 + \tilde{\gamma}_2 <1$ and in particular $\beta_1, \gamma_1, \beta_2, \gamma_2, \tilde{\beta}_1, \tilde{\gamma}_1, \tilde{\beta}_2, \tilde{\gamma}_2 \in (0,1)$, we can treat the four integral products and the remaining eight addenda of the right-hand side in a such way that eventually they are absorbed by the three integral terms involving the gradients in the left one. More exactly, to the products we apply the first inequality derived in Lemma \ref{LemmaIneq}, and to the other terms Young's inequality. In this way, the resulting linear combination of 
the terms $\int_{\Omega} |\nabla |\nabla v|^q|^2$, $\int_{\Omega} |\nabla |\nabla w|^r|^2$ and $\int_{\Omega} |\nabla (u+1)^{\frac{m_1+p-1}{2}}|^2$ can be written as 
$\frac{c_{54}}{2} \int_{\Omega} |\nabla |\nabla v|^q|^2 + \frac{c_{55}}{2} \int_{\Omega} |\nabla |\nabla w|^r|^2 
+\frac{c_{56}}{2} \int_{\Omega} |\nabla (u+1)^{\frac{m_1+p-1}{2}}|^2$, which throughout relation \eqref{Somma1} infers the claim.
\end{proof}
\end{lemma}

\begin{remark}\label{Limitaz}
We observe that the argument of Lemma \ref{Met_GN} can be applied to the linear case $m_1=m_2=m_3=1$ only for $\alpha \in \left(0,\frac{2}{n}\right)$ and/or 
$\gamma \in \left(0,\frac{2}{n}\right)$.
\end{remark}

\subsection{The logistic case}\label{Log}
For the logistic case we retrace part of the computations above connected to the usage of the Young inequality only.
\begin{lemma}\label{Estim_general_For_u^p_nablav^2qLemmaLog}
If $m_1,m_2,m_3 \in \R$ comply with $m_1>\max\left\{2m_2-1, 2m_3-1, \frac{n-2}{n}\right\}$ or  $m_1>\max\left\{m_2-\frac{1}{n}, m_3-\frac{1}{n}, \frac{n-2}{n}\right\}$ or $m_1>\max\left\{2m_2-1, m_3-\frac{1}{n}, \frac{n-2}{n}\right\}$ or $m_1>\max\left\{m_2-\frac{1}{n}, 2m_3-1, \frac{n-2}{n}\right\}$ or $m_1>\max\left\{2m_2-\beta, 2m_3-\beta, \frac{n-2}{n}\right\}$ or $m_1>\max\{2m_2-\beta, 2m_3-\beta\}$ whenever $\alpha, \gamma \in \left(0,\frac{1}{n}\right]$, or $m_1>\max\left\{2m_2, 2m_3, \frac{n-2}{n}\right\}$ 
or $m_1>\max\{2m_2+1-\beta, 2m_3+1-\beta\}$ whenever $\alpha, \gamma \in \left(\frac{1}{n},1\right)$, then there exist $p, q, r>1$ such that $(u,v,w)$ satisfies 
a similar inequality as in \eqref{MainInequality}.
\begin{proof}
As in Lemma \ref{Estim_general_For_u^p_nablav^2qLemma}, in view of inequalities \eqref{Y0} and \eqref{Y1} taking into account \eqref{Young} and \eqref{Young1}, and the properties of the logistic $h$ in \eqref{h}, relation \eqref{Estim_1_For_u^p} now becomes for some positive $\tilde{c}_{3}$ and for all $t \in (0,\TM)$
\begin{equation}\label{Estim_1_For_u^p1}  
\begin{split}
\frac{d}{dt} \int_\Omega (u+1)^p & \leq (-p(p-1)+\delta_1) \int_\Omega (u+1)^{p+m_1-3} |\nabla u|^2 \\
&\quad\,+ \tilde{c}_1 \int_\Omega (u+1)^{\frac{(p+2m_2-m_1-1)s}{s-2}}+ \tilde{c}_2 \int_\Omega (u+1)^{\frac{(p+2m_3-m_1-1)s}{s-2}}\\
&\quad\,+ pk_+ \int_\Omega (u+1)^p - p \mu  \int_\Omega (u+1)^{p-1}u^{\beta} + \tilde{c}_{3}. 
\end{split}
\end{equation}
Applying the inequality $(A+B)^p \leq 2^{p-1} (A^p+B^p)$ with $A,B \geq 0$ and $p>1$ to the last integral in \eqref{Estim_1_For_u^p1}, implies that  
$-u^{\beta} \leq -\frac{1}{2^{\beta-1}} (u+1)^{\beta}+1$; therefore 
\begin{equation}\label{beta}
- p \mu  \int_\Omega (u+1)^{p-1} u^{\beta} \leq -\frac{p \mu}{2^{\beta-1}} \int_\Omega (u+1)^{p-1+\beta} + p \mu \int_\Omega (u+1)^{p-1} \quad \textrm{on } (0,\TM).
\end{equation}
Henceforth, by taking into account the Young inequality, we have that for $t \in (0,\TM)$
\begin{align} \label{k}
\begin{split}
&pk_+ \int_\Omega (u+1)^p \leq \delta_1 \int_\Omega (u+1)^{p-1+\beta} + \tilde{c}_{4} \quad \textrm{and} \\
&p \mu \int_\Omega (u+1)^{p-1} \leq \delta_2 \int_\Omega (u+1)^{p-1+\beta} + \tilde{c}_{5},
\end{split}
\end{align}
with $\delta_1, \delta_2 >0$ and some $\tilde{c}_4, \tilde{c}_5 >0$.\\
\textbf{Case 1}: $\alpha, \gamma \in (0, \frac{1}{n}]$ and $m_1>\max\{2m_2-1, 2m_3-1, \frac{n-2}{n}\}$ or  $m_1>\max\{m_2-\frac{1}{n}, m_3-\frac{1}{n}, \frac{n-2}{n}\}$ or $m_1>\max\{2m_2-1, m_3-\frac{1}{n}, \frac{n-2}{n}\}$ or $m_1>\max\{m_2-\frac{1}{n}, 2m_3-1, \frac{n-2}{n}\}$ or $m_1>\max\{2m_2-\beta, 2m_3-\beta, \frac{n-2}{n}\}$ or $m_1>\max\{2m_2-\beta, 2m_3-\beta\}$. 
For $m_1>\max\left\{2m_2-1, 2m_3-1, \frac{n-2}{n}\right\}$ or  $m_1>\max\left\{m_2-\frac{1}{n}, m_3-\frac{1}{n}, \frac{n-2}{n}\right\}$ or $m_1>\max\left\{2m_2-1, m_3-\frac{1}{n}, \frac{n-2}{n}\right\}$ or $m_1>\max\left\{m_2-\frac{1}{n}, 2m_3-1, \frac{n-2}{n}\right\}$, we refer to Lemma \ref{Estim_general_For_u^p_nablav^2qLemma} and we take in mind inequalities \eqref{Young4}, \eqref{GN2}, \eqref{GNY}, \eqref{Young1}, \eqref{Yxi} and \eqref{GNY1}.
Conversely, when $m_1>2 m_2-\beta$ and $m_1>2 m_3-\beta$, we have that (recall $s$ may be arbitrary large) $\frac{(p+2m_2-m_1-1)s}{s-2} < p-1+\beta$ 
and $\frac{(p+2m_3-m_1-1)s}{s-2} < p-1+\beta$, and by means of the Young inequality estimates \eqref{Young4} and \eqref{Yxi} can alternatively read
\begin{equation}\label{young4}  
\tilde{c}_1 \int_\Omega (u+1)^{\frac{(p+2m_2-m_1-1)s}{s-2}} \leq \delta_3 \int_\Omega (u+1)^{p-1+\beta} + \tilde{c}_6 \quad \textrm{for all } t \in(0,\TM),
\end{equation}
and 
\begin{equation}\label{young5}  
\tilde{c}_2 \int_\Omega (u+1)^{\frac{(p+2m_3-m_1-1)s}{s-2}} \leq \delta_4 \int_\Omega (u+1)^{p-1+\beta} + \tilde{c}_7 \quad \textrm{on }\, (0,\TM),
\end{equation}
with $\delta_3, \delta_4>0$ and positive $\tilde{c}_6, \tilde{c}_7$.
By inserting estimates \eqref{beta} and \eqref{k} into relation \eqref{Estim_1_For_u^p1}, as well as taking into account \eqref{Young4} and \eqref{Yxi}
(or, alternatively to \eqref{Young4} and \eqref{Yxi}, bound \eqref{young4} and \eqref{young5}), for suitable $\hat{\epsilon}, \tilde{\delta} >0$ and some $\tilde{c}_7>0$ we arrive at 
\begin{align*}
\frac{d}{dt} \int_\Omega (u+1)^p &\leq \left(-\frac{4p(p-1)}{(m_1+p-1)^2} + \hat{\epsilon}\right) \int_\Omega |\nabla (u+1)^{\frac{m_1+p-1}{2}}|^2 \\
&\quad\,+ \left(\tilde{\delta} - \frac{p \mu}{2^{\beta-1}} \right) \int_\Omega (u+1)^{p-1+\beta} + \tilde{c}_7 \quad \textrm{for all } t \in (0,\TM),
\end{align*}
where we used again relation \eqref{GradU}. We can conclude reasoning exactly as in the second part of the proof of Lemma \ref{Estim_general_For_u^p_nablav^2qLemma},
by exploiting $m_1>\frac{n-2}{n}$ and by choosing suitable $\hat{\epsilon},\tilde{\delta}$. On the other hand, by enlarging $p$, Young's inequality allows us to obtain the following alternative estimates
\begin{align*}
c_{21} \int_\Omega u^{2\alpha} |\nabla v|^{2q-2} &\leq \delta_5 \int_\Omega u^{p-1+\beta} + \tilde{c}_8 \int_\Omega |\nabla v|^{\frac{2(q-1)(p-1+\beta)}{p-1+\beta-2\alpha}}\\
&\leq \delta_5 \int_\Omega u^{p-1+\beta} + \delta_6 \int_\Omega |\nabla v|^s + \tilde{c}_9 \quad \textrm{on } (0,\TM),
\end{align*}
and 
\begin{align*}
c_{27} \int_\Omega u^{2\gamma} |\nabla w|^{2r-2} &\leq \delta_6 \int_\Omega u^{p-1+\beta} + \tilde{c}_9 \int_\Omega |\nabla w|^{\frac{2(r-1)(p-1+\beta)}{p-1+\beta-2\gamma}}\\
&\leq \delta_6 \int_\Omega u^{p-1+\beta} + \delta_7 \int_\Omega |\nabla w|^s + \tilde{c}_{10} \quad \textrm{for all } t \in (0,\TM).
\end{align*}
\textbf{Case 2}: $\alpha, \gamma \in \left(\frac{1}{n},1\right)$ and $m_1>\max\left\{2m_2, 2m_3, \frac{n-2}{n}\right\}$ 
or $m_1>\max\{2m_2+1-\beta, 2m_3+1-\beta\}$.\\
For $m_1>\max\left\{2m_2, 2m_3, \frac{n-2}{n}\right\}$ we will refer to the second case of Lemma \ref{Estim_general_For_u^p_nablav^2qLemma}. 
Now, if $m_1>2m_2+1-\beta$ and $m_1>2m_3+1-\beta$, then some $p$ sufficiently large infers to $\frac{(p+2m_2-m_1-1)(p+1)}{p}< p-1+\beta$ and 
$\frac{(p+2m_3-m_1-1)(p+1)}{p}< p-1+\beta$, so that for any positive $\delta_8, \delta_9$ 
and some $\tilde{c}_{10}, \tilde{c}_{11}>0$ we have
\[
\bar{c}_1 \int_\Omega (u+1)^{\frac{(p+2m_2-m_1-1)(p+1)}{p}} \leq \delta_8 \int_\Omega (u+1)^{p-1+\beta} + \tilde{c}_{10} \quad \textrm{on } (0,\TM),
\]
and 
\[
\bar{c}_4 \int_\Omega (u+1)^{\frac{(p+2m_3-m_1-1)(p+1)}{p}} \leq \delta_9 \int_\Omega (u+1)^{p-1+\beta} + \tilde{c}_{11} \quad \textrm{for all } t \in(0,\TM).
\]
Now, the integrals $\int_\Omega u^{2\alpha} |\nabla v|^{2p-2}$ and $\int_\Omega u^{2\gamma} |\nabla w|^{2p-2}$ can be treated as in 
\textbf{Case~2} of Lemma \ref{Estim_general_For_u^p_nablav^2qLemma} or alternatively, by exploiting $\alpha, \gamma <1$, a different application of Young's inequalities leads on $(0,\TM)$ to
\begin{equation}\label{Pannocchia}
	\begin{split}
c_{21} \int_\Omega u^{2\alpha} |\nabla v|^{2p-2} &\leq \delta_5 \int_\Omega u^{p-1+\beta} + \tilde{c}_8 \int_\Omega |\nabla v|^{\frac{2(p-1)(p-1+\beta)}{p-1+\beta-2\alpha}}\\
&\leq \delta_5 \int_\Omega u^{p-1+\beta} + \delta_{10} \int_\Omega |\nabla v|^{2(p+1)} + \tilde{c}_{12},
\end{split}
\end{equation}
and 
\begin{equation}\label{Disumana}
	\begin{split}
c_{27} \int_\Omega u^{2\gamma} |\nabla w|^{2p-2} &\leq \delta_6 \int_\Omega u^{p-1+\beta} + \tilde{c}_9 \int_\Omega |\nabla w|^{\frac{2(p-1)(p-1+\beta)}{p-1+\beta-2\gamma}}\\
&\leq \delta_6 \int_\Omega u^{p-1+\beta} + \delta_{11} \int_\Omega |\nabla w|^{2(p+1)} + \tilde{c}_{13},
\end{split}
\end{equation}
with $\delta_{10}, \delta_{11}>0$ and some positive $\tilde{c}_{12}, \tilde{c}_{13}$.
The remaining part of the proof follows as \textbf{Case 2} of Lemma \ref{Estim_general_For_u^p_nablav^2qLemma} for the terms dealing with 
$\int_\Omega |\nabla v|^{2(p+1)}$ and $\int_\Omega |\nabla w|^{2(p+1)}$.

As before, this result applies also for
\begin{itemize}
	\setlength{\itemsep}{0.7mm}
	\item [$\rhd$] $\alpha \in \left(0,\frac{1}{n}\right], \gamma \in \left(\frac{1}{n},1\right)$ and $m_1>\max\Bigl\{2m_2-1,\frac{n-2}{n}, 2m_3\Bigr\}$ or 
	$m_1>\max\Bigl\{2m_2-1,\frac{n-2}{n}, 2m_3+1-\beta\Bigr\}$ or $m_1>\max\Bigl\{m_2-\frac{1}{n},2m_3,\frac{n-2}{n}\Bigr\}$ or
	$m_1>\max\Bigl\{m_2-\frac{1}{n},2m_3+1-\beta\Bigr\}$ or $m_1>\max\Bigl\{2m_2-\beta,2m_3,\frac{n-2}{n}\Bigr\}$ or
	$m_1>\max\Bigl\{2m_2-\beta,2m_3+1-\beta,\frac{n-2}{n}\Bigr\}$ or $m_1>\max\Bigl\{2m_2-\beta,2m_3+1-\beta\Bigr\}$,
\item [$\rhd$] $\alpha \in \left(\frac{1}{n},1\right), \gamma \in \left(0,\frac{1}{n}\right]$ and 
	$m_1>\max\Bigl\{2m_2, 2m_3-1,\frac{n-2}{n}\Bigr\}$ or $m_1>\max\Bigl\{2m_2+1-\beta, 2m_3-1,\frac{n-2}{n}, \Bigr\}$ 
	or $m_1>\max\Bigl\{2m_2, m_3-\frac{1}{n},\frac{n-2}{n}\Bigr\}$ or $m_1>\max\Bigl\{2m_2+1-\beta, m_3-\frac{1}{n}\Bigr\}$
	or $m_1>\max\Bigl\{2m_2, 2m_3-\beta,\frac{n-2}{n}\Bigr\}$ or $m_1>\max\Bigl\{2m_2+1-\beta, 2m_3-\beta,\frac{n-2}{n}\Bigr\}$ 
	or $m_1>\max\Bigl\{2m_2+1-\beta,2m_3-\beta\Bigr\}$.
\end{itemize}
\end{proof}
\end{lemma}
As a by-product of what now obtained we are in a position to conclude.
\subsection{Proof of Theorems \ref{MainTheorem} and \ref{MainTheorem1}}
\begin{proof}
Let $(u_0,v_0,w_0) \in (W^{1,\infty}(\Omega))^3$ with $u_0, v_0, w_0 \geq 0$ on $\bar{\Omega}$. For $f$ and $g$ as in \eqref{f} and, respectively, for $f$, $g$ as in \eqref{f} and $h$ as in \eqref{h}, let $\alpha, \gamma >0$ and let $m_1,m_2,m_3 \in \R$ comply with \ref{A1}---\ref{A16}, respectively, \textcolor{blue}{\ref{A17}---\ref{A20}}. Then, we refer to Lemmas \ref{Estim_general_For_u^p_nablav^2qLemma} and \ref{Met_GN}, respectively, Lemma \ref{Estim_general_For_u^p_nablav^2qLemmaLog} and obtain for some $C_1,C_2,C_3, C_4>0$
\begin{align}\label{Estim_general_For_y_2}
\begin{split}
&y'(t) + C_1 \int_\Omega |\nabla (u+1)^{\frac{m_1+p-1}{2}}|^2 \\
&\quad\,+ C_2 \int_\Omega |\nabla |\nabla v|^q|^2 + C_3 \int_\Omega |\nabla |\nabla w|^r|^2\leq C_4
\quad \textrm{ on } (0, \TM).
\end{split}
\end{align}
Successively, the Gagliardo--Nirenberg inequality again makes that some positive constants $c_{62}, c_{63}, c_{64}$ imply from the one hand
\begin{equation}\label{G_N2}
\int_\Omega (u+1)^p \leq c_{62} \Big(\int_\Omega |\nabla (u+1)^\frac{m_1+p-1}{2}|^2\Big)^{\kappa_1}+ c_{62}  \quad \textrm{for all } t \in (0,\TM),
\end{equation}
(as already done in the inequality immediately before \eqref{GN2}), and from the other on $(0,\TM)$
\begin{align}\label{Estim_Nabla nabla v^q} 
\begin{split}
\int_\Omega \lvert \nabla v\rvert^{2q}
&=\lvert \lvert \lvert \nabla v\rvert^q\lvert \lvert_{L^2(\Omega)}^2\\
&\leq c_{63} \lvert \lvert\nabla  \lvert \nabla v \rvert^q\rvert \lvert_{L^2(\Omega)}^{2\kappa_2} \lvert \lvert\lvert \nabla v \rvert^q\lvert \lvert_{L^\frac{1}{q}(\Omega)}^{2(1-\kappa_2)}+c_{63} \lvert \lvert \lvert \nabla v \rvert^q\lvert \lvert^2_{L^\frac{1}{q}(\Omega)},
\end{split}
\end{align}
and similarly for all $t \in (0,\TM)$
\begin{align}\label{Estim_Nabla nabla w^r}
\begin{split}
\int_\Omega \lvert \nabla w\rvert^{2r}
&=\lvert \lvert \lvert \nabla w\rvert^r\lvert \lvert_{L^2(\Omega)}^2 \\
&\leq c_{64} \lvert \lvert\nabla  \lvert \nabla w \rvert^r\rvert \lvert_{L^2(\Omega)}^{2\kappa_3} \lvert \lvert\lvert \nabla w \rvert^r\lvert \lvert_{L^\frac{1}{r}(\Omega)}^{2(1-\kappa_3)}+c_{64} \lvert \lvert \lvert \nabla w \rvert^r\lvert \lvert^2_{L^\frac{1}{r}(\Omega)}, 
\end{split}
\end{align}
with $\kappa_2, \kappa_3$ already defined in Lemma \ref{LemmaCoefficientAiAndExponents}. 
Subsequently, the $L^s$-bound of $\nabla v$ in \eqref{Cg} and of $\nabla w$ in \eqref{Cg1} infer some $c_{65}, c_{66}>0$ such that
\begin{equation*}\label{Estim_Nabla nabla v^p^2}
\int_\Omega \lvert \nabla v\rvert^{2q}\leq c_{65} \Big(\int_\Omega \lvert \nabla \lvert \nabla v \rvert^q\rvert^2\Big)^{\kappa_2}+c_{65} \quad \textrm{on } (0,\TM),
\end{equation*}
and
\begin{equation*}\label{Estim_Nabla nabla w^p^2}
\int_\Omega \lvert \nabla w\rvert^{2r}\leq c_{66} \Big(\int_\Omega \lvert \nabla \lvert \nabla w \rvert^r\rvert^2\Big)^{\kappa_3}+c_{66} \quad \textrm{for all } t \in (0,\TM).
\end{equation*}
At this stage, by using  estimates \eqref{G_N2}, \eqref{Estim_Nabla nabla v^q} and \eqref{Estim_Nabla nabla w^r}, with the aid of the second inequality in Lemma \ref{LemmaIneq}, relation \eqref{Estim_general_For_y_2} entails positive constants  $c_{67}$ and $c_{68}$, and $\tilde{\kappa}=\min\{\frac{1}{\kappa_1},\frac{1}{\kappa_2}, \frac{1}{\kappa_3}\}$ such that 
\begin{equation*}\label{MainInitialProblemWithM}
\begin{cases}
y'(t)\leq c_{67}-c_{68} y^{\tilde{\kappa}}(t)\quad \textrm{on } (0,\TM),\\
y(0)=\int_\Omega (u_0+1)^p+ \int_\Omega |\nabla v_0|^{2q} + \int_\Omega |\nabla w_0|^{2r}. 
\end{cases}
\end{equation*}
Finally, ODE comparison principles imply $u \in L^\infty((0,\TM);L^p(\Omega))$, and the conclusion is a consequence of the boundedness criterion in \eqref{BoundednessCriterio}.
\end{proof}

\begin{remark}\label{DettagliCasoLineareNonLogistico}
First we note that the conclusions \ref{NonLogistic2} and \ref{NonLogistic3} in Remark \ref{CommentiNonLogistico} can be similarly justified by using reasonings given above and connected to the proof of Theorem \ref{MainTheorem}. In particular, two issues are crucial: the regularity of $\nabla v, \nabla w$ (see Lemma \ref{LocalV}) and the problematic integral term $\int_\Omega u^{p+1}$, which can be treated by applying the Gagliardo--Nirenberg and Young's inequalities for $n=1$ without any condition, and for $n=2$ provided some restrictions. 

As to the other cases, we can discuss these scenarios:
\begin{itemize}
\setlength{\itemsep}{0.7mm}
\item If $\alpha,\gamma \in \left[\frac{2}{n},1\right)$ (\ref{A3}) we cannot apply Lemma \ref{Met_GN} (see Remark \ref{Limitaz}). Moreover, even 
Lemma \ref{Estim_general_For_u^p_nablav^2qLemma} does not work; in fact, from it we would obtain that $\alpha, \gamma \in [2,1)$ for $n=1$ and 
$\alpha, \gamma \in [1,1)$ for $n=2$, which is not possible. The same contradictions appear for the cases $\alpha\in\left[\frac{2}{n},1\right]$, $\gamma \in \left[\frac{2}{n},1\right)$ (\ref{A4}) or $\alpha \in \left[\frac{2}{n},1\right)$, $\gamma \in \left[\frac{2}{n},1\right]$ (\ref{A5}) and $\alpha,\gamma \in \left[\frac{2}{n},1\right]$ (\ref{A6}). 
\item If $\alpha\in\left(0,\frac{1}{n}\right]$, $\gamma \in \left[\frac{2}{n},1\right)$ (\ref{A8}) or $\alpha \in \left[\frac{2}{n},1\right)$, $\gamma \in \left(0,\frac{1}{n}\right]$ (\ref{A13}), even though both Lemmas \ref{Estim_general_For_u^p_nablav^2qLemma} and \ref{Met_GN} are applicable, a further contradiction appears. 
\item If $\alpha\in\left(\frac{1}{n},\frac{2}{n}\right)$, $\gamma \in \left[\frac{2}{n},1\right)$ (\ref{A11}) or $\alpha \in \left[\frac{2}{n},1\right)$, 
$\gamma \in \left(\frac{1}{n},\frac{2}{n}\right)$ (\ref{A15}), and similarly for  $\alpha\in\left(\frac{1}{n},\frac{2}{n}\right)$, $\gamma \in \left[\frac{2}{n},1\right]$ (\ref{A12}) or $\alpha \in \left[\frac{2}{n},1\right]$, 
$\gamma \in \left(\frac{1}{n},\frac{2}{n}\right)$ (\ref{A16}),  by reasoning as in the previous cases we obtain a contradiction, characterized by the fact that the intervals of $\alpha$ or $\gamma$ are empty.
\end{itemize}
\end{remark}
\begin{remark}\label{DettagliCasoLineareLogistico}
	Let us spend some words on the hints mentioned in Remark \ref{CommentiLogistico}.
	\begin{itemize}
	\item [$\triangleright$] For the linear case $m_1=m_2=m_3=1$, a simple substitution ensures the validity of Theorem \ref{MainTheorem1} for $\beta>2$ and $\alpha,\gamma\in (0,1)$; conversely,  for  $\alpha,\gamma\in (0,1]$ relations \eqref{Pannocchia} and \eqref{Disumana} are still applicable for $\beta>2$ (item \ref{logistico1}); for the nonlinear case, the same reasoning can be carried out to show \ref{logistico5}, \ref{logistico6} and \ref{logistico7};
		\item [$\triangleright$] For $\beta=2$, in bound \eqref{Estim_1_For_u^p1}, the term associated to the logistic dampening effect takes the form $-p\mu \int_\Omega u^{p+1}$, so that for $\mu$ large, the other positive contributions proportional as well to  $\int_\Omega u^{p+1}$ itself, can be absorbed. As to the expressions of $\mu$, relations in \ref{logistico2}, \ref{logistico3} and \ref{logistico4} come from the related range of $\alpha$ and $\gamma$; for $\beta=2$ and the nonlinear case, a further largeness assumption on $\mu$ is required (\ref{logistico8}, \ref{logistico9} and \ref{logistico10}). 
\end{itemize}
\end{remark}


\begin{thebibliography}{10}

\bibitem{BaghaeiKhelghatiIranian}
K.~Baghaei and A.~Khelghati.
\newblock Boundedness of classical solutions for a chemotaxis model with
  consumption of chemoattractant.
\newblock {\em C. R. Math. Acad. Sci. Paris}, 355(6):633--639, 2017.

\bibitem{ChiyoMarrasTanakaYokota2021}
Y.~Chiyo, M.~Marras, Y.~Tanaka, and T.~Yokota.
\newblock Blow-up phenomena in a parabolic-elliptic-elliptic
  attraction-repulsion chemotaxis system with superlinear logistic degradation.
\newblock {\em Nonlinear Anal.}, 212:Paper No. 112550, 14, 2021.

\bibitem{FrassuGalvanViglialoro}
S.~Frassu, R.~{Rodr\'{i}guez Galv\'{a}n}, and G.~Viglialoro.
\newblock Uniform in time ${L}^{\infty}$-estimates for an attraction-repulsion
  chemotaxis system with double saturation.
\newblock {\em Discrete Contin. Dyn. Syst. - B.}, \url{ arXiv:2106.07038v2},
  2022.

\bibitem{FrassuCorViglialoro}
S.~Frassu, C.~{van der Mee}, and G.~Viglialoro.
\newblock Boundedness in a nonlinear attraction-repulsion {K}eller--{S}egel
  system with production and consumption.
\newblock {\em J. Math. Anal. Appl.}, 504(2):125428, 2021.

\bibitem{frassuviglialoro1}
S.~Frassu and G.~Viglialoro.
\newblock Boundedness for a fully parabolic {K}eller--{S}egel model with
  sublinear segregation and superlinear aggregation.
\newblock {\em Acta Appl. Math.}, 171(1):1--20, 2021.

\bibitem{frassuviglialoroConsumptionProduction}
S.~Frassu and G.~Viglialoro.
\newblock Boundedness in a chemotaxis system with consumed chemoattractant and
  produced chemorepellent.
\newblock {\em Nonlinear Anal.}, 213:112505, 2021.

\bibitem{GuoJiangZhengAttr-Rep}
Q.~Guo, Z.~Jiang, and S.~Zheng.
\newblock Critical mass for an attraction-repulsion chemotaxis system.
\newblock {\em Appl. Anal.}, 97(13):2349--2354, 2018.

\bibitem{K-S-1970}
E.~F. Keller and L.~A. Segel.
\newblock Initiation of slime mold aggregation viewed as an instability.
\newblock {\em J. Theoret. Biol.}, 26(3):399--415, 1970.

\bibitem{Keller-1971-MC}
E.~F. Keller and L.~A. Segel.
\newblock {Model for chemotaxis.}
\newblock {\em J. Theoret. Biol.}, 30(2):225--234, 1971.

\bibitem{Keller-1971-TBC}
E.~F. Keller and L.~A. Segel.
\newblock {Traveling bands of chemotactic bacteria: A theoretical analysis}.
\newblock {\em J. Theoret. Biol.}, 30(2):235, 1971.

\bibitem{LankeitWangConsumptLogistic}
J.~Lankeit and Y.~Wang.
\newblock Global existence, boundedness and stabilization in a high-dimensional
  chemotaxis system with consumption.
\newblock {\em Discrete Contin. Dyn. Syst.}, 37(12):6099--6121, 2017.

\bibitem{LI-LiAttrRepuls}
Y.~Li and Y.~Li.
\newblock Blow-up of nonradial solutions to attraction-repulsion chemotaxis
  system in two dimensions.
\newblock {\em Nonlinear Anal. Real World Appl.}, 30:170--183, 2016.

\bibitem{Luca2003Alzheimer}
M.~Luca, A.~Chavez-Ross, L.~Edelstein-Keshet, and A.~Mogilner.
\newblock {Chemotactic signaling, microglia, and Alzheimer's disease senile
  plaques: Is there a connection?}
\newblock {\em Bull. Math. Biol.}, 65(4):693--730, 2003.

\bibitem{MarrasViglialoroMathNach}
M.~Marras and G.~Viglialoro.
\newblock Boundedness in a fully parabolic chemotaxis-consumption system with
  nonlinear diffusion and sensitivity, and logistic source.
\newblock {\em Math. Nachr.}, 291(14--15):2318--2333, 2018.

\bibitem{TaoBoun}
Y.~Tao.
\newblock Boundedness in a chemotaxis model with oxygen consumption by
  bacteria.
\newblock {\em J. Math. Anal. Appl.}, 381(2):521--529, 2011.

\bibitem{TaoWanM3ASAttrRep}
Y.~Tao and Z.-A. Wang.
\newblock Competing effects of attraction vs. repulsion in chemotaxis.
\newblock {\em Math. Models Methods Appl. Sci.}, 23(1):1--36, 2013.

\bibitem{TaoWinkParaPara}
Y.~Tao and M.~Winkler.
\newblock {Boundedness in a quasilinear parabolic-parabolic Keller--Segel
  system with subcritical sensitivity}.
\newblock {\em J. Differerential Equations}, 252(1):692--715, 2012.

\bibitem{VIGLIALORO-JMAA-BlowUp-Attr-Rep}
G.~Viglialoro.
\newblock Explicit lower bound of blow-up time for an attraction-repulsion
  chemotaxis system.
\newblock {\em J. Math. Anal. App.}, 479(1):1069--1077, 2019.

\bibitem{ViglialoroMatNacAttr-Repul}
G.~Viglialoro.
\newblock Influence of nonlinear production on the global solvability of an
  attraction-repulsion chemotaxis system.
\newblock {\em Math. Nachr.}, 294(12):2441--2454, 2021.

\bibitem{WinklerN-Sto_CPDE}
M.~Winkler.
\newblock {Global large-data solutions in a chemotaxis-(Navier--)Stokes system
  modeling cellular swimming in fluid drops}.
\newblock {\em Comm. Partial Differential Equations}, 37(2):319--351, 2012.

\bibitem{WinklerN-Sto_2d}
M.~Winkler.
\newblock { Stabilization in a two-dimensional chemotaxis-Navier--Stokes
  system}.
\newblock {\em Arch. Ration. Mech. Anal.}, 212(2):455--487, 2014.

\end{thebibliography}
\medskip

\end{document}